\documentclass[12pt]{article}
\usepackage{geometry}                
\usepackage{graphicx,color,mathtools}
\usepackage{amssymb,amsmath,amsthm,mathrsfs}
\usepackage[all,cmtip]{xy}
\usepackage{epstopdf, comment}
\usepackage{bm,enumitem} 

\numberwithin{equation}{section}

\newtheorem{theorem}{Theorem}[section]

\newtheorem{lemma}[theorem]{Lemma}
\newtheorem{corollary}[theorem]{Corollary}

\theoremstyle{remark}
\newtheorem*{remark}{Remark}
\newtheorem*{example}{Example}

\theoremstyle{definition}

\DeclareMathOperator{\supp}{supp }
\DeclareMathOperator{\diam}{diam}
\DeclareMathOperator{\aut}{Aut}

\DeclareMathOperator{\hyp}{hyp}
\DeclareMathOperator{\dist}{dist}
\DeclareMathOperator{\area}{Area}
\DeclareMathOperator{\PSL}{PSL}
\DeclareMathOperator{\loc}{loc}

\DeclareMathOperator{\Mod}{Mod}
\DeclareMathOperator{\Root}{root}
\DeclareMathOperator{\parent}{parent}

\DeclareMathOperator{\ord}{ord}

\DeclareMathOperator{\Haus}{H}

\DeclareMathOperator{\Right}{right}
\DeclareMathOperator{\Left}{left}

\DeclareMathOperator{\twig}{tw}
\DeclareMathOperator{\rad}{rad}

\DeclareMathOperator{\Ann}{\bf A}
\DeclareMathOperator{\Rect}{\bf R}
\DeclareMathOperator{\pRect}{\partial \bf R}
\DeclareMathOperator{\Lip}{Lip}

\DeclareMathOperator{\re}{Re}
\DeclareMathOperator{\im}{Im}

\DeclareMathOperator{\length}{length}

\DeclareMathOperator{\hypTriangle}{\triangle_{\text{hyp}}}
\DeclareMathOperator{\sph}{\hat{\mathbb{C}}}
 
\newcommand{\D}{\mathbb{D}}
\newcommand{\T}{\mathcal T}
\newcommand{\CC}{\mathbb{C}}

\newcommand{\TPL}{\rm{TPL}_1}

\newcommand{\TW}{\mathscr{TW}}
\newcommand{\QS}{\rm{QS}}

\title{Shapes of infinite conformally balanced trees}
\author{Oleg Ivrii,\, Peter Lin,\, Steffen Rohde,\, Emanuel Sygal}

\date{October 31, 2023}

\begin{document}

\maketitle

\begin{abstract}
 Numerical experiments by Werness, Lee and the third author suggested that dessin d'enfants associated to large trivalent trees approximate the developed deltoid introduced by Lee, Lyubich, Makarov and Mukherjee. In this paper, we confirm this conjecture. As a side product of our techniques, we give a new proof of a theorem of Bishop which says that ``true trees are dense.'' We also exhibit a sequence of trees whose conformally natural shapes converge to the cauliflower, the Julia set of $z\mapsto z^2+1/4$.
\end{abstract}

\section{Introduction}

A finite tree $\mathcal T$ in the plane is called a {\em conformally balanced tree} or a {\em true tree} if 

\begin{enumerate}[leftmargin=1.75cm, label={\rm (TT\arabic*)}]
\item \label{item:true-tree1} Every edge has the same harmonic measure as seen from infinity.

\item \label{item:true-tree2} Harmonic measures on the two sides of every edge are identical.
\end{enumerate}
 Conformally balanced trees are in one-to-one correspondence with Shabat polynomials: any conformally balanced tree is the pre-image of the segment $[-1,1]$ by an essentially unique  polynomial $p$ with critical values $\pm 1$. (The polynomial $p(z)$ is uniquely determined up to multiplication by $-1$.)
  
We say that two trees $T_1, T_2$ in the plane are  {\em equivalent} if there is an orientation-preserving homeomorphism of the plane which takes $T_1$ onto $T_2$. It is a classic fact that every finite tree $T$ in the plane is equivalent to a conformally balanced tree $\mathcal T$, which is unique up to post-composition with affine maps. A proof of these facts will be sketched in Section \ref{sec:shabat}.

It is natural to ask if infinite trees also have a natural shape. In \cite{lin-rohde}, the second- and third author developed the theory of Gehring trees and showed that the Aldous continuum random tree  possesses a natural conformal structure.
In this paper, we consider the {\em infinite trivalent tree}\/ $\mathcal T$, which exhibits a different and surprising behaviour. To come up with a natural shape for $\mathcal T$, we truncate it at level $n$, form the conformally balanced tree $\mathcal T_n$ and take $n \to \infty$.

In order for the finite trees $\mathcal T_n$ to converge, we need to normalize them in some way. Throughout the rest of the paper, we use the {\em hydrodynamic} normalization: we ask that each conformal map $\varphi_n: \sph \setminus \mathbb{D} \to \sph \setminus \mathcal T_n$ has the expansion $z \to z + O(1/z)$ near infinity.

Our main theorem states:

\begin{theorem}
  \label{main-thm}
  The trees $\mathcal T_n$ converge in the Hausdorff topology to an infinite trivalent tree union a Jordan curve
  $\mathcal T_\infty \cup \partial \Omega$. The domain $\Omega$ enclosed by $\partial \Omega$ is the developed deltoid.  The Shabat polynomials $p_n$ converge to $F \circ R^{-1}$ where $F$ is a modular function invariant under an index 2 subgroup of $\PSL(2, \mathbb{Z})$ and $R: \mathbb{D} \to \Omega$ is the Riemann map.
\end{theorem}

\begin{figure}[h]
\centering
   \includegraphics[scale=0.18]{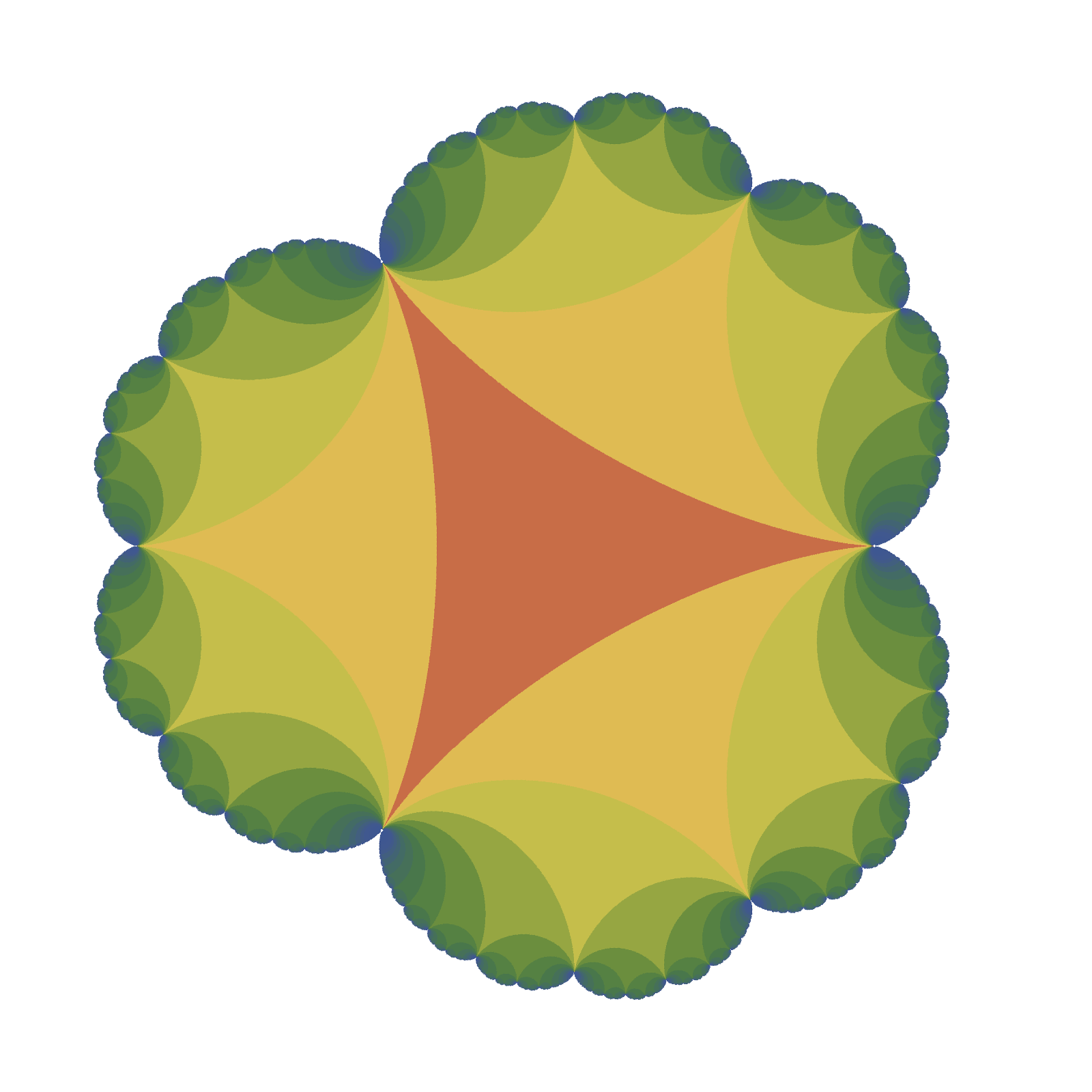} \quad \includegraphics[scale=0.18]{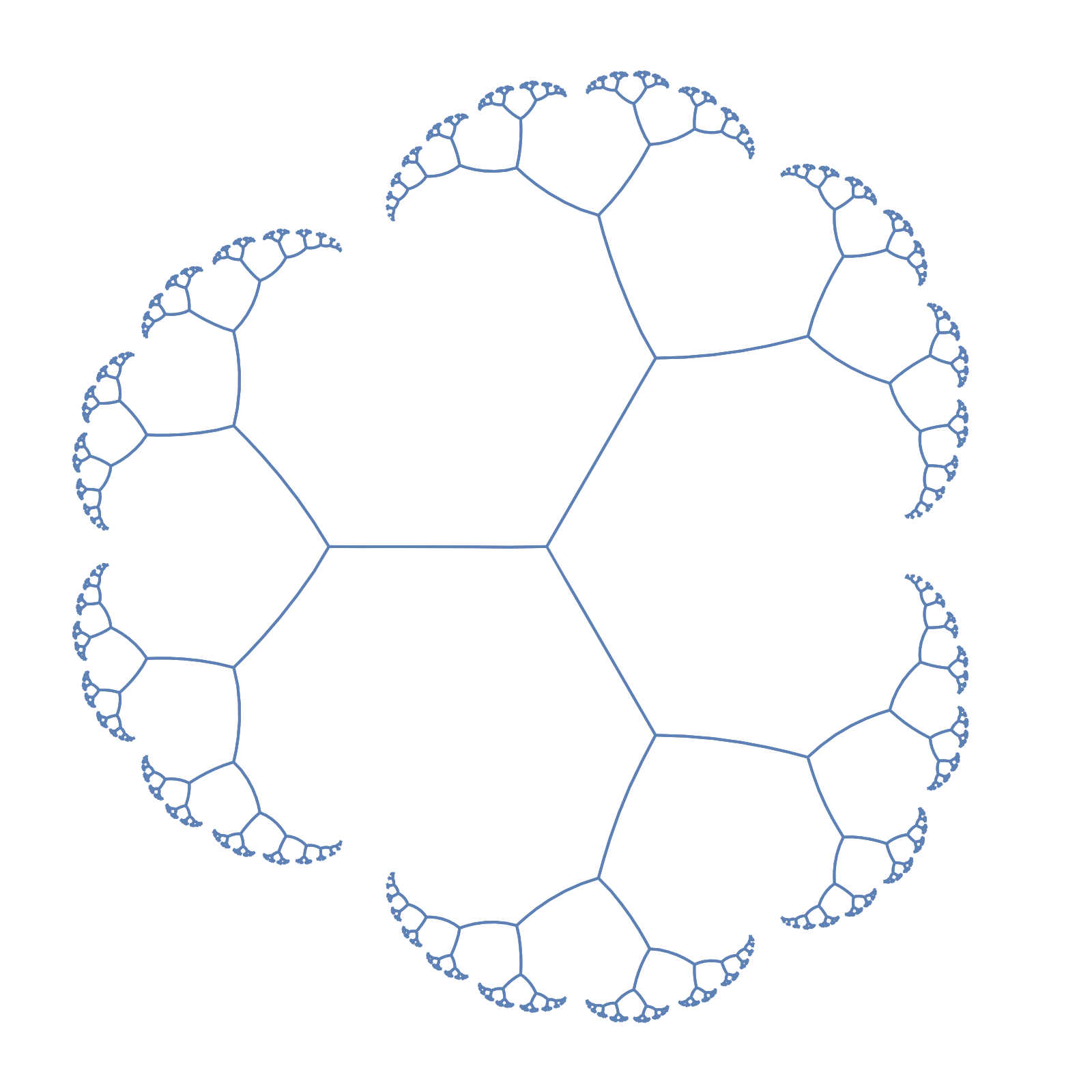}
 \caption{The developed deltoid and its approximating conformally balanced tree.}
 \label{fig:developed-deltoid}
\end{figure}

\begin{remark}
The choice of truncation is important: by considering other truncations of the infinite trivalent tree, one can obtain different limit sets. In fact, any compact connected set in the plane can be approximated in the Hausdorff topology by conformally balancing finite truncations of the infinite trivalent tree, thereby giving another proof of a theorem of Bishop \cite{bishop}. See Appendix \ref{sec:true-trees-are-dense}.
\end{remark}

\subsection{The developed deltoid}\label{subsec:deltoid}

The {\em deltoid} $\triangle \subset \mathbb{C}$ is a remarkable domain in the plane bounded by a Jordan curve with three outward pointing cusps. It can be described as the curve traversed by a point on a circle of radius $1/3$ as it rolls around in the interior of a circle of radius 1. 
Alternatively, one can describe the exterior of the deltoid $\triangle_e = \mathbb{C} \setminus \overline{\triangle}$ as the image of $\mathbb{D}_e = \sph \setminus \mathbb{D}$ under the conformal map
$z \to z + \frac{1}{2z^2}$.

The exterior of the deltoid is part of a somewhat mysterious family of domains called {\em quadrature domains}.  Quadrature domains have several equivalent definitions such as possessing a {\em Schwarz reflection} which is an anti-holomorphic function $\sigma: \triangle_e \to \mathbb{C}$ that is identity on $\partial \triangle_e$. By repeatedly reflecting the deltoid in its sides, one obtains the {\em developed deltoid} 
$$
	\Omega = \bigcup_{k \ge 0} \sigma^{-k}(\triangle),
$$
see Fig.~\ref{fig:developed-deltoid}. The developed deltoid was first studied by S-Y.~Lee, M.~Lyubich, N.~G.~Makarov and S.~Mukherjee \cite{schwarz-reflections}, who showed that it fuses Fuchsian dynamics with anti-holomorphic dynamics:

\begin{theorem}
{\em (i)} The boundary of the developed deltoid $\partial \Omega$ is the unique Jordan curve that realizes the mating of the ideal triangle group and $z \to \overline{z}^2$.

{\em (ii)} The developed deltoid $\Omega$ is a John domain. In particular, $\partial \Omega$ is conformally removable.
\end{theorem}

We now describe their result in detail:

\medskip

	{\em Dynamics on $\mathbb{D}_e$.} In the exterior of the unit disk, we consider the dynamical system $z \to \overline{z}^2$.

\medskip

	{\em Dynamics on $\mathbb{D}$.}  Let $\triangle_{\hyp} \subset \mathbb{D}$ be the ideal triangle in the unit disk with vertices at $1, \omega, \omega^2$, where $\omega = e^{2\pi i/3}$ is a third root of unity.  Consider the group $\Gamma = \langle R_{\rho_1}, R_{\rho_2}, R_{\rho_3} \rangle \subset \aut(\mathbb{D})$ generated by the reflections in the sides $s_1, s_2, s_3$ of $\triangle_{\hyp}$. The images
$$
\{\gamma(\triangle_{\hyp}) : \gamma \in \Gamma\}
$$ tessellate the unit disk.
The Markov map $\rho: \mathbb{D} \setminus \triangle_{\hyp} \to \mathbb{D}$ is defined as $R_{\rho_1}$ on the (hyperbolic) half-plane cut off by $s_1$, $R_{\rho_2}$ on the half-plane cut off by $s_2$ and $R_{\rho_3}$ on the half-plane cut off by $s_3$.

\medskip

	{\em What it means to be a mating.} A Jordan curve $\gamma = \partial \Omega$ is a {\em mating} of $z \to \overline{z}^2$ and $\Gamma$ if there exist conformal maps $\varphi: \mathbb{D} \to \Omega$, $\psi: \mathbb{D}_e \to \Omega_e$ that glue the dynamical systems together, i.e.~$\varphi \circ \rho \circ \varphi^{-1} = \psi \circ \overline{z}^2 \circ \psi^{-1}$ on $\partial \Omega$. In particular, this implies that
$$
  \sigma(z) =
  \begin{cases}
    \psi \circ \overline{z}^2 \circ \psi^{-1}, \quad z \in \overline{\Omega_e} \\
    \varphi \circ \rho \circ \varphi^{-1}, \quad z \in \overline{\Omega} \setminus \varphi(\triangle_{\hyp})
  \end{cases}
$$
is a  Schwarz reflection for $\sph  \setminus \varphi(\triangle_{\hyp})$, and hence $\sph  \setminus \varphi(\triangle_{\hyp})$ is a quadrature domain.

\medskip

A set $E$ is called {\em conformally removable} if every conformal map $h: \sph \setminus E \to \sph \setminus F$ which extends continuously to the Riemann sphere is a M\"obius transformation.

\subsection{Strategy of proof}

Our proof of Theorem \ref{main-thm} proceeds in three steps:

\medskip

{\em Step 1.} We first show that any subsequential limit of the true trees $\mathcal T_n$ in the Hausdorff topology is homeomorphic to an infinite trivalent tree union a Jordan curve $\mathcal T_\infty \cup \partial \Omega$, with $\mathcal T_\infty \subset \Omega$. Among our key tools are estimates for the diameters of edges by means of conformal modulus estimates of certain curve families. A notable difference to the setting of random trees is that in the truncated trivalent tree, the diameters of a fixed edge do not shrink to zero as $n\to\infty$, see also the remark at the end of Section \ref{sec:interior-edges}.

\medskip

{\em Step 2.} We then show that any subsequential limit $\partial \Omega$ realizes the mating of $z \to \overline{z}^2$ and $\Gamma$. At this point, one can appeal to the uniqueness of the mating \cite{schwarz-reflections} to complete the proof of Theorem \ref{main-thm}. However, appealing to  \cite{schwarz-reflections} feels somewhat unsatisfactory since it relies on a priori knowledge of the deltoid, while ideally, one would want to ``discover'' the deltoid from the infinite trivalent tree. 

\medskip

{\em Step 3.}
To show that the limit of the $\mathcal T_n$ does not depend on the subsequence, we prove ``partial conformal removability.'' Partial conformal removability is a much less stringent property than full conformal removability and it is easier to check. In essence, it asks that if $h: \sph \setminus E \to \sph \setminus F$ is a conformal map (which extends continuously to the Riemann sphere) onto the complement of a set $F$ which has roughly the same geometry as $E$, then $h$ is a M\"obius transformation.

\subsection{Acknowledgements}
In 2014, Brent Werness (oral communication) proposed to study the natural shape of the infinite trivalent tree and posed the question of the ``shrinking of diameters.'' During a visit of Seung-Yeop Lee to Seattle in 2015, Brent, Seung-Yeop and the third author performed computer experiments and observed the similarity between the trees and the developed deltoid, leading to the conjecture regarding their convergence. We are grateful to Brent and Seung-Yeop for our discussions and their contributions. We would also like to thank Curt McMullen for his continued interest in this project and for his suggestion that Shabat polynomials converge to a modular function.

This research was supported by the Israeli Science Foundation (3134/21) and the National Science Foundation (DMS-1700069 and DMS-1954674).

\section{Preliminaries}

In this section, we gather a number of useful facts that will be used in this paper. We also describe the Farey tessellation and discuss weak conformal removability.

\subsection{Moduli of annuli and rectangles}

It is well known that any doubly-connected domain $\Ann \subset \mathbb{C}$ can be mapped onto a round annulus $\{ z : r < |z| < R\}$. The number $\Mod \Ann :=\frac{1}{2\pi} \log \frac{R}{r}$ is called the {\em modulus} of $\Ann$. Two doubly-connected domains are conformally equivalent if and only if their moduli coincide.

 A {\em metric} $\rho(z)$ is a non-negative measurable function defined on a domain $\Omega \subset \mathbb{C}$. One can use $\rho(z)$ to measure lengths of rectifiable curves
$$
\ell_{\rho}(\gamma) = \int_\gamma \rho(z) |dz|
$$
and compute areas of shapes, for instance the total area of $\rho$ is given by
$$
A(\rho) = \int_{\Omega} \rho(z)^2 |dz|^2.
$$

The metric $\rho$ is said to be {\em admissible} for a family of rectifiable curves $\Gamma$ contained in $\Omega$ if the $\rho$-length of every curve $\gamma \in \Gamma$ is at least 1. The {\em modulus} of the curve family $\Gamma$ is defined as
$$
\Mod \Gamma := \inf_\rho A(\rho),
$$
where the infimum is taken over all admissible metrics $\rho$. If one finds a conformal metric $\rho$ such that $\ell_{\rho}(\gamma) \ge L$ for any $\gamma \in \Gamma$, then $\Mod \Gamma \le A(\rho)/L^2$.

The modulus of a doubly-connected domain is a special case of the above construction: $\Mod \Ann$ is equal to the modulus of the family of curves $\Gamma_\circlearrowleft$ that separate the two boundary components, while
$1/\Mod \Ann$ is equal to modulus of the family $\Gamma_{\uparrow}$ of curves  that connect the opposite boundary components of $\Ann$. Thus one uses $\Gamma_{\circlearrowleft}$ to give upper bounds for $\Mod \Ann$ while one uses $\Gamma_{\uparrow}$ go give lower bounds for $\Mod \Ann$.

We will frequently use the following two simple rules for modulus, which follow from the definitions:

\begin{enumerate}

\item (Monotonicity rule) If $\Ann_1 \subset \Ann$ is an essential doubly-connected subdomain, so that $\Gamma_\circlearrowleft(\Ann_1) \subset \Gamma_\circlearrowleft(\Ann)$, then
$\Mod \Ann_1 \le \Mod \Ann$.

\item (Parallel rule) If a doubly-connected domain $\Ann = \Ann_1 \cup \Ann_2$ can be represented as a union of two essential doubly-connected domains, then 
$$\Mod \Ann_1 + \Mod \Ann_2 \le \Mod \Ann\!.$$
\end{enumerate}

We will also use the following standard estimates: 

\begin{lemma}
\label{lem:egg_yolk}

 Let $\Omega $ be a simply-connected domain in the plane.

{\em (a)} Suppose $F$ is a compact connected set contained in $\Omega$. If $\Mod (\Omega \setminus F) \ge m$ is bounded from below, then $$\dist(\partial \Omega, F) \ge c \, \diam F,$$
for some $c > 0$ which depends only on $m > 0$.  Furthermore, $c\to\infty$ as $m\to\infty.$
Conversely, if $\dist(\partial \Omega, F) \ge c \, \diam F,$ then $\Mod (\Omega \setminus F) \ge m(c).$

{\em (b)} Suppose $E \subset F$ are two compact connected sets contained in $\Omega$. If 
$$m_1 \, \le \, \Mod (\Omega \setminus E) \, \le \, \Mod (\Omega \setminus F) \, \le \, m_2, $$ then $\diam E \asymp \diam F.$ 
In fact, there exists a constant $C = C(m_1, m_2) > 1$ so that
$
F \subset B(e,\, C \cdot \diam E)
$
for any point $e \in E$, where $B(x,r)$ denotes the ball of radius $r$ centered at $x.$
\end{lemma}

A {\em conformal rectangle} $\Rect$ is a simply connected domain with four marked prime ends $z_1, z_2, z_3, z_4$. In this paper, all conformal rectangles will be {\em marked}\/, i.e.~equipped with a  distinguished pair of opposite sides.
The Schwarz-Cristoffel formula provides a conformal map from $\Rect$  onto a geometric rectangle $[0,m] \times [0,1]$. If one insists that the marked sides of $\Rect$ are mapped onto the vertical sides of $[0,m] \times [0,1]$, then the number $m \in (0,\infty)$ is determined uniquely. The number $m := \Mod \Rect$ is known as the {\em modulus} of $\Rect$ and is
equal to the modulus of the curve family $\Gamma_{\updownarrow}$ which separates the distinguished pair of opposite sides.

For further properties of conformal modulus, we refer the reader to \cite[Chapter 4]{garnett-marshall} and \cite[Chapter 2]{mcmullen}.

\subsection{Farey tesellation}
\label{sub:Farey}

\begin{figure}[h]
\centering
   \includegraphics[scale=.3]{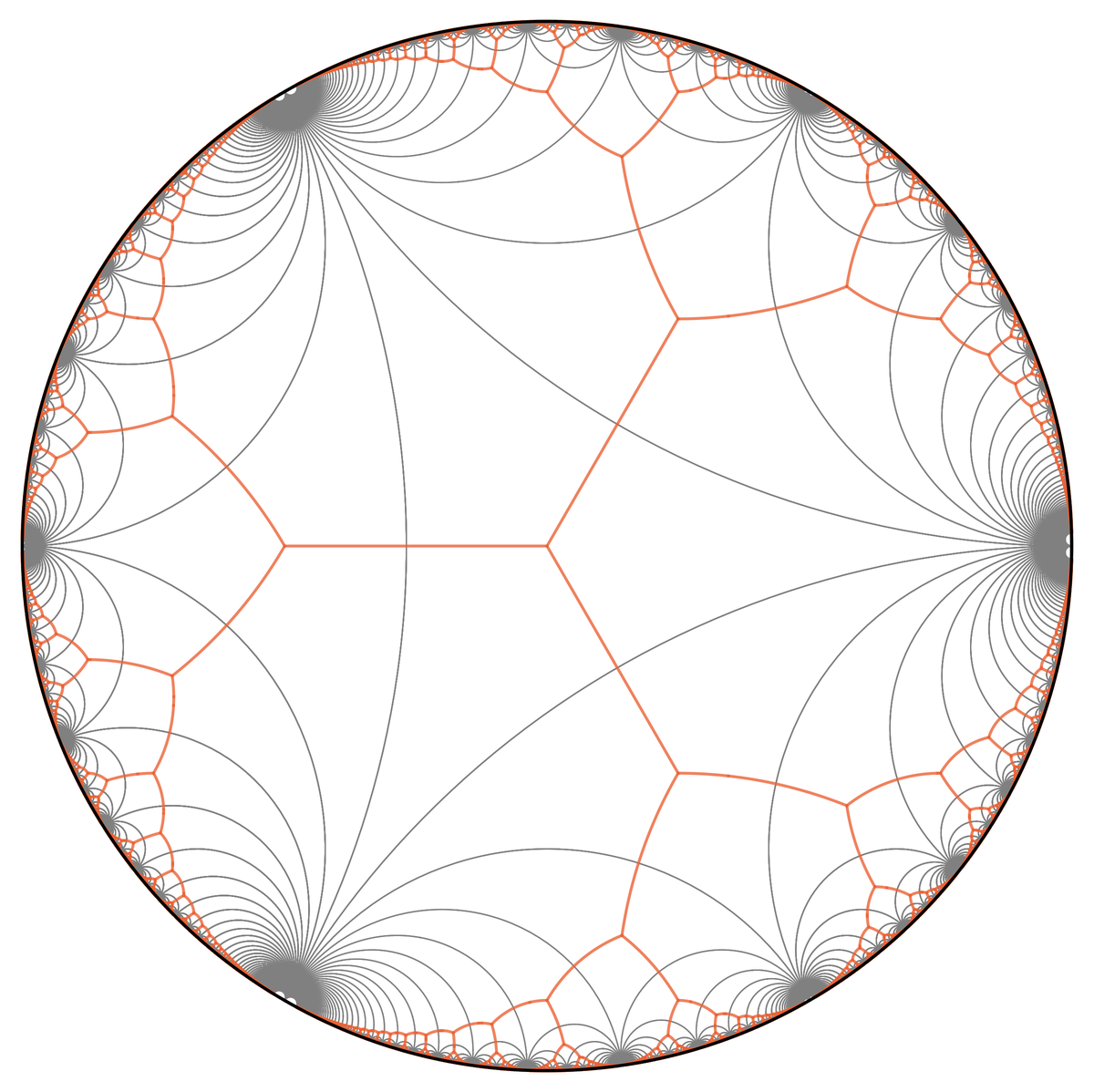}
   \caption{The Farey tesselation and Farey tree}
    \label{fig:farey}
\end{figure}

Let $\hypTriangle \subset \mathbb{D}$ be the ideal triangle in the unit disk with vertices
$1$, $\omega = e^{2\pi i /3}$ and $\overline{\omega} = e^{4\pi i /3}$. Repeatedly reflecting $\hypTriangle$ in its sides, one obtains a tessellation of the unit disk by ideal triangles. The dual graph (which joins centers of the triangles by hyperbolic geodesics) is called the {\em Farey tree} $\mathcal F$, see Figure \ref{fig:farey}.

We designate the center of $\hypTriangle$ as the root vertex. Each non-root triangle $\triangle$ can be labeled by a digit $1,2,3$ followed by a finite sequence of $L$'s and $R$'s, which indicates the path one travels from $\hypTriangle$ to $\triangle$.
For example, in the word
$$
2 \, \underbrace{L}_{k_1 = 1} \, \underbrace{R}_{k_2 = 1} \,  \underbrace{L L L}_{k_3=3} \, \underbrace{R R}_{k_4=2} \, \underbrace{L L}_{k_5=2} \, \underbrace{R}_{k_6=1} \, \underbrace{L L L L L}_{k_7=5} \, \underbrace{R R}_{k_8=2},
$$
 the digit 2 indicates that we start by walking along the dual tree from the root vertex to its second child. After the first step, each vertex has two children and we have to decide whether to turn left or right. The options are indicated by `$L$' and `$R$' respectively.

\begin{lemma}
For a non-root triangle $\triangle$  in the Farey tessellation,
$$
\log \frac{1}{\diam \triangle} \asymp \sum_{i=1}^m  \log(1+k_i).
$$
\end{lemma}

\begin{proof}
It is easier and clearly equivalent to work in the upper half plane $\mathbb{H}$ where $\hypTriangle$ has vertices 0, 1 and $\infty$ and in the first step, we walk down. 
Let 
$$\triangle_0 = \hypTriangle, \ \triangle_1 = (0, 1/2, 1), \ \triangle_2, \  \dots, \ \triangle_n = \triangle$$ be the sequence of triangles from $\hypTriangle$ to $\triangle$.
Each triangle $\triangle_j$ in this sequence has three vertices on the real axis $a_j < b_j < c_j$. To estimate $\diam \triangle_j$,
we keep track of the ratio $$r(\triangle_j) := \frac{b_j - a_j }{ c_j - a_j},
$$
which measures the distortion of the triangle $\triangle_j$.
Each time we do an right turn after a left turn or vice versa, the ratio is ``reset'' to a value in $[1/3, 2/3]$. After a series of $k$ consecutive left turns,  $1 - r \asymp 1/k$, while after a series of $k$ consecutive right turns, $r \asymp 1/k$.

After making $k$ left or right turns in a row, the diameter goes down by a factor of roughly $k+1$: for $1 \le k \le k_{j+1}$,
$$
\log \frac{1}{\diam \triangle_{k_1 +k_2 + \dots + k_j + k}} -\log \frac{1}{ \diam \triangle_{k_1+k_2+\dots+k_j+1}}  \, \asymp \, \log(k+1).
$$
  When we make a right turn after a series of $k_j$ left turns (or a left turn after a series of $k_j$ left turns), the diameter goes down by a factor of $k_j+1$, i.e.~
$$
\log \frac{1}{\diam \triangle_{k_1 +k_2 + \dots + k_j + 1}} -\log \frac{1}{ \diam \triangle_{k_1+k_2+\dots+k_j}}  \,  \asymp\, \log (k_j+1).
$$
The above equations give the desired bound for $\diam \triangle$.
\end{proof}

\subsection{Weak conformal removability}

Suppose $X$ and $X'$ are two compact sets in the complex plane and $\varphi: \sph \setminus X \to \sph \setminus X'$ is a conformal map that extends continuously to a homeomorphism of the sphere.
We describe a condition which guarantees that $\varphi$ is a M\"obius transformation:

\begin{lemma}
\label{weak-conformal-removability}
Suppose that there is a countable exceptional set $E \subset X$ and a countable collection of closed subsets $s_1, s_2, \dots$ of $X$, called shadows, such that every point in $X \setminus E$ belongs to infinitely many sets $s_i$. If
\begin{equation}
\label{eq:square-sum-shadows}
\sum_{i=1}^\infty \diam^2 s_i < \infty, \qquad \sum_{i=1}^\infty \diam^2 \varphi(s_i) < \infty,
\end{equation}
then $\varphi$ is a M\"obius transformation.
\end{lemma}

 For convenience, we write $s'_i = \varphi(s_i)$. Note that (\ref{eq:square-sum-shadows}) implies that $X$ and $X'$ have 2-dimensional Lebesgue measure 0.

\begin{proof}
Call a direction $v$ {\em good} if for almost every line $\mathcal \ell$ pointing in the direction of $v$, the set $\varphi(\ell \cap X)$ has linear Lebesgue measure 0. One says that $\varphi$ is {\em absolutely continuous on lines} (ACL)  if the  directions parallel to the coordinate axes are good. It is well known that if $\varphi \in W_{\loc}^{1,2}(\mathbb{C} \setminus X)$ is ACL, then $\varphi \in W_{\loc}^{1,2}(\mathbb{C})$. Weyl's lemma then guarantees that $\varphi$ is conformal on the Riemann sphere, and therefore, a M\"obius transformation. Below, we will show that every direction is good. 

Instead of showing that a set has zero 1-dimensional Lebesgue measure $m_1$, we may instead show that it has zero 1-dimensional content $m_1^\infty$. The definition of 1-dimensional content is similar to that of 1-dimensional measure, but allows covers by balls of arbitrary size. Therefore, the lemma reduces to showing that for almost every line $\ell$ parallel to $v$,  the 1-dimensional content of $\varphi(\ell \cap X)$ is 0.

Since the set $E$ is countable, almost every line $\ell$ parallel to $v$ misses $E$. For such a line,
\begin{equation}
\label{eq:desired-estimate}
m_1^\infty(\varphi(\ell \cap X)) \le \sum_{s_i \cap \ell \ne \emptyset, \,i > N} 
\diam s'_i.
\end{equation}
The last equation holds for {\em any} $N \ge 1$ since any point in $X \setminus E$ is contained in infinitely many shadows, which allows us to avoid putting the first $N-1$ shadows in the cover.
In other words,
\begin{equation}
\label{eq:shadow-sum}
 \sum_{s_i \cap \ell \ne \emptyset} 
\diam s'_i < \infty \quad \implies \quad  m_1^\infty(\varphi(\ell \cap X)) = 0.
\end{equation}
As
\begin{align*}
\int_{\ell || v} \biggl \{ \sum_{s_i \cap \ell \ne \emptyset} 
\diam s'_i \biggr \} \, d\ell & \, \le \, \sum_{i=1}^\infty 
\diam s_i \cdot \diam s'_i \\
& \, \le \, \frac{1}{2} \biggl ( \sum_{i=1}^\infty 
\diam^2 s_i + \diam^2 s'_i \biggr ) \\
& \, < \, \infty,
\end{align*}
 the integrand must be finite for a.e.~$\ell$. This completes the proof.
\end{proof}

\section{Background on true trees}

In this section, we discuss the link between true trees and Shabat polynomials. We then describe the local geometry of true trees whose vertices have bounded valence. Finally, we define shortcuts and obstacles that will be used to give moduli estimates to control the global geometry of trees.

\subsection{True trees and Shabat polynomials}
\label{sec:shabat}

Let $T$ be a finite tree in the plane.
To find its conformally balanced shape $\mathcal T$, label the sides of  edges of $T$ in counter-clockwise order: $\vec{e}_1, \vec{e}_2, \dots, \vec{e}_{2N}$.
For each half-edge $\vec{e}_i$, form an equilateral triangle $\triangle(\vec{e}_i, \infty)$ whose sides have unit length.

We first glue these equilateral triangles into a $2N$-gon ${\bm D}_{2N}$ with sides $\vec{e}_i$, labeled counter-clockwise, and central vertex $\infty$. We then glue $\vec{e}_i$ with $\vec{e}_j$ whenever $\vec{e}_i, \vec{e}_j$ are opposite sides of the same edge $e \in T$. This construction produces a topological sphere which has a flat structure away from the cone points at the
vertices of the triangles. Uniformizing this sphere produces the desired tree $\mathcal T \subset \hat{\mathbb{C}}$.

Associated to a true tree $\mathcal T$ is a {\em Shabat polynomial} $p(z)$ with critical values $\pm 1$ such that $\mathcal T = p^{-1}([-1,1])$.
To construct $p$, colour each triangle $\triangle(\vec{e}_i, \infty) \subset \hat{\mathbb{C}} \setminus \mathcal T$ black or white, so that adjacent triangles have opposite colours.  On each
black triangle $\triangle(\vec{e}_i, \infty)$, define $p(z)$ to be the conformal map onto the upper half-plane $\mathbb{H}$ which takes $\vec{e}_i \to [-1,1]$ and $\infty \to \infty$. Similarly, on each white triangle $\triangle(\vec{e}_i, \infty)$,
define $p$ to be the conformal map onto the lower half-plane $\mathbb{L}$ which takes $\vec{e}_i \to [-1,1]$ and $\infty \to \infty$.

Since $\mathcal T$ is a true tree, $p$ extends to a continuous function on the Riemann sphere.
As $\mathcal T$ is made up of real-analytic arcs, $p$ is meromorphic on the Riemann sphere, and hence a rational function. As the only pole of $p$ is at infinity, it is
a polynomial. Finally, since $p$ is $N:1$ at infinity, $p$ is a polynomial of degree $N$. From the construction, it is readily seen that $p$ has critical values $\pm 1$
and $\mathcal T = p^{-1}([-1,1])$.

In order to define the Shabat polynomial uniquely, we need to specify which vertices are sent to $+1$ and which vertices are sent to $-1$. Making a different choice amounts to multiplying $p(z)$ by $-1$.
If $\mathcal T$ has a distinguished vertex $v_{\Root}$, then it is natural to choose the Shabat polynomial so that
$p(v_{\Root}) = 1$.

\subsection{Trees of bounded valence}
\label{sec:trees-bounded-valence}

We now present some general results on the local behaviour of true trees. The following lemmas say  that true trees whose vertices have bounded valence are well-behaved: neighbouring edges have comparable size and
the relative distance between non-adjacent edges is bounded below.

\begin{lemma}
  \label{basic-trees1}
  Let $d \ge 2$ be an integer. Suppose $e = \overline{v_1v_2}$ is an edge in a true tree $\mathcal T$ with $\deg v_1 \le d$ and $\deg v_2 \le d$.
  There is a simply connected neighbourhood $U \supset e$ with $\Mod (U \setminus e) \ge m(d)$
  such that only edges adjacent to $e$ can intersect $U$.
\end{lemma}

\begin{lemma}
  \label{basic-trees2}
  Fix an integer $d \ge 2$. Let $v$ be a vertex of a conformally balanced tree $\mathcal T$. If the degrees of all vertices in
  $
    \{ w:  d_{\mathcal T}(v,w) \le 2 \}
  $
  are $\le d$, then the diameters of the edges $\overline{vv_i}$ emanating from $v$ are comparable (with the comparison constant depending on $d$).
\end{lemma}

The proofs use the concept of a {\em star} of a vertex in a true tree. For a vertex $v$ of $\mathcal T$, we define $\star_v$ as the union of the triangles $\triangle(\vec{e}, \infty)$ that contain $v$. We enumerate
the $2 \deg v$ triangles  in $\star_v$
counter-clockwise: $\triangle_1, \triangle_2, \dots, \triangle_{2 \deg v}$.

Now, decompose the unit disk $\mathbb{D}$ into $2\deg v$ sectors $\sigma_1, \sigma_2, \dots, \sigma_{2\deg v}$ using $2 \deg v$ equally-spaced radial rays.
For each $i = 1, 2, \dots, 2 \deg v$, let $\psi_i$ be the conformal map from $\sigma_i$ to $\triangle_i$ which takes vertices to vertices, with $0$ mapping to $v$.
Since $\mathcal T$ is a true tree, the maps $\psi_i$ glue along radial rays to form a conformal map
$
  \psi_v: \mathbb{D} \to \star_v.
$

On an edge $e_i = \overline{vv_i}$ of $\mathcal T$ emanating from $v$, we mark the points $a_i, b_i$ such that the segments $\overline{va_i}, \overline{a_i b_i}, \overline{b_iv_i}$ have equal length in the equilateral triangle model of $\triangle(\vec{e}_i, \infty)$. Note that the points $a_i, b_i$ do not depend on which one of the two sides of $e_i$ is used.

Applying Koebe's distortion theorem to $\psi_v$ tells us that the diameters of the $2\deg v$ segments
$$
  \bigl \{ \overline{va_i}, \ \overline{a_i b_i} \, : \, i = 1, 2, \dots, \deg(v) \bigr \}
$$
are comparable. By considering stars centered at the neighbouring vertices $v_i$, we see that
$$
  \bigl \{ \overline{a_i b_i}, \ \overline{b_iv_i} \, : \,  i = 1, 2, \dots, \deg(v) \bigr \}.
$$
are also comparable. Putting these estimates together proves Lemma \ref{basic-trees2}.

Lemma \ref{basic-trees1} follows from Lemma \ref{lem:egg_yolk} (a) after applying Koebe's distortion theorem to $\psi_{v_1}$ and $\psi_{v_2}$. Similar reasoning shows:

\begin{lemma}
  Suppose $\{\mathcal T_n\}_{n=0}^\infty$ is an infinite sequence of conformally balanced trees whose vertices have uniformly bounded degrees. Then any subsequential Hausdorff limit of a sequence of edges $e^{(n)} \subset \mathcal T_n$ is either a point or a real-analytic arc.
\end{lemma}

\begin{proof}
  Suppose the edge $e^{(n)}$ connects the vertices $v_1^{(n)}$ and $v_2^{(n)}$.  As above, we mark the points $a^{(n)}$ and $b^{(n)}$ which trisect the edge $e^{(n)}$. We pass to a subsequence so that the maps $\psi_{v_1}^{(n)}$ and $\psi_{v_2}^{(n)}$ converge uniformly on compact subsets of the unit disk.

  If the limiting maps $\psi_{v_1} = \lim_{n \to \infty} \psi_{v_1}^{(n)}$ and $\psi_{v_2} = \lim_{n \to \infty}  \psi_{v_2}^{(n)}$ are constant, then the edges $e^{(n)}$ collapse to a point. Otherwise, the limiting edge $e = \lim e^{(n)}$ is covered by two compatible real-analytic arcs $\overline{v_1 b} = \lim_{n \to \infty} \overline{v_1^{(n)} b^{(n)}}$ and $\overline{a v_2} = \lim_{n \to \infty} \overline{ a^{(n)} v_2^{(n)}}$.
\end{proof}

\subsection{Shortcuts and obstacles}
\label{sec:shortcuts-and-ostacles}

Let $\mathcal T$ be a conformally balanced tree in the plane, normalized so that the Riemann map
$\varphi: \hat{\mathbb{C}} \setminus \overline{\mathbb{D}} \to \hat{\mathbb{C}} \setminus \mathcal T$ satisfies
$\varphi(z) = z + O(1/z)$ as $z \to \infty$.

To control the geometry of $\mathcal T$, we estimate conformal moduli of various path families $\Gamma$ contained in doubly-connected domains $\Ann \subset\CC$.
 An instructive example is the family of closed curves surrounding an edge of the tree, which will be discussed in detail in Section \ref{sec:interior-edges}.

The idea behind our estimates is as follows: since we will only estimate moduli in the setting of finite balanced trees, we will not have to worry about the possibility that the area of $\T$ might be positive. By conformal invariance, we may estimate the modulus in any of the three conformally equivalent models $\sph\setminus \mathcal T$,  ${\bm D}_{2N}/\!\!\sim$ or 
$(\sph\setminus \D)/\!\!\sim$. In the latter model, the equivalence relation on $\partial\D$ is given by the identifications of $\varphi$ and the family $\varphi^{-1}(\Gamma)$ consists of sets $\varphi^{-1}(\gamma)$ that may be disconnected: if a curve $\gamma\in\Gamma$ crosses an edge $e\in\T,$ then $\varphi^{-1}(\gamma)$ {\em enters} one side of $\varphi^{-1}(e)$,  {\em teleports}  through the identification provided by $\varphi$, and {\em exits} on the other side of $\varphi^{-1}(e)$.

We will construct admissible metrics of the form
\begin{equation}
 \label{eq:metric}
\rho = \alpha_0 \Bigl( \rho_0 + \sum_{e\in\T}  \alpha_e  \rho_e \Bigr),
\end{equation}
where the {\em background metric}  $\rho_0={\bf1}_{\varphi^{-1}(\Ann)}$ serves the purpose of controlling the length of curves $\gamma$ that do not intersect $\T$, while the
{\em obstacles} $\rho_e$
have the purpose of penalizing teleportation so that shortcuts are not worthwhile. The constant $\alpha_0$ is chosen so that curves that do not intersect $\T$ have length $\geq1$ under $ \alpha_0 \rho_0.$ 

We build the obstacles $\rho_e$ so that they assign length $\geq1$ to all curves $\gamma$ that intersect $e$ (and are not confined to the union of the triangles that are incident to $e$). It is easiest to describe the construction in ${\bm D}_{2N}/\!\sim$, which is a surface composed of $2N$ equilateral triangles $\triangle(\vec{e}_i, \infty)$ of side length 1: namely, we define
$\rho_e$  as  three times the characteristic function of the 1/3-neighborhood of $e$ in the flat metric, i.e.~
$$\rho_e =  3\times  {\bf1}_{B_{1/3}(e)},$$
where $B_{1/3}(e)$ is the set of points of distance at most $1/3$ from $e$.

We denote the conformal transport of this metric to $\hat{\mathbb{C}} \setminus \overline{\mathbb{D}}$
again by $\rho_e$. Since any point $z \in \hat{\mathbb{C}} \setminus \overline{\mathbb{D}}$ can be in the support of at most $D = \max_{v\in\T} \deg(v)$ obstacles, it can be in at most $D+1$ of the sets $\supp \rho_0 \cup \{\supp \rho_{e}\}$, and the area of $\rho$ can be estimated by
\begin{equation}
  \label{eq:area}
A(\rho) \, \le \, (D+1)^2\, \alpha_0^2 \Bigl(  A(\rho_0) + \sum_e \alpha_e^2 \cdot A(\rho_e)^2 \Bigr) \, \lesssim \,
\alpha_0^2 \Bigl( A(\rho_0) + \sum_e \alpha_e^2 \Bigr).
\end{equation}

For an edge $e$ in $\mathcal T$, we denote by $\mathcal T(e) \subset \mathcal T$ the subtree consisting of the edge $e$ and its descendants (as measured from the root vertex).
It is easy to see that $S(e) = \varphi^{-1}(\mathcal T(e))$ is an arc in the unit circle $\partial\D$. We define the {\em outer shortcut} of $e$ as the Euclidean length of $S(e)$:
$$
s(e) = \length (S(e)).
$$
We define  $\mathcal T^{-}(e) = \mathcal T(e) \setminus e$ as the union of all the descendants of $e$. Naturally, we define $S^{-}(e) = \varphi^{-1}(\mathcal T^{-}(e))$ and $s^{-}(e) =  \length (S^{-}(e)) = s(e) - \pi/N$. Unless $e$ is a boundary edge, the difference between the outer and inner shortcuts is not significant.

\subsection{A lower bound for the diameters of edges}
\label{sec:interior-edges}

Let $\mathcal T$ be a true tree and
$\Omega_2 \subset \mathbb{C}$ be the simply-connected domain bounded by the equipotential curve
$
  \varphi \bigl ( \{ z: |z| = 2\} \bigr ).
$
The hydrodynamic normalization of the conformal map $\varphi$ implies that $\diam \mathcal T \geq c_0 > 0$ is bounded from below by a universal constant (the sharp value $c_0=2$ is irrelevant for our purpose).

In view of Lemma \ref{lem:egg_yolk}, to give a lower bound for the diameter of an edge $e_0$ in $\mathcal T$, it is enough to give an upper bound for the modulus of the family of curves $\Gamma_{\circlearrowleft}(\Ann)$ that separate the boundary components of $\Ann = \Omega_2 \setminus e_0$. We will now show that the metric
$$
\rho = \frac{1}{s^{-}(e_0)} \Bigl( {\bf 1}_{A(0;1,2)} + \sum_{e\in \mathcal T^{-}_n(e_0)} s(e)\rho_e\Bigr)
$$ 
is admissible for $\varphi^{-1}(\Gamma_{\circlearrowleft}(\Ann)),$ where the summation is over the descendants of $e_0.$
Consider a curve $\gamma\in\Gamma_{\circlearrowleft}(\Ann)$. If we pull $\gamma$ back by $\varphi^{-1}$, we  get a path in the annulus
$$A(0;1,2) = \{z : 1 \le |z| < 2 \},$$ which may teleport from $x \in  \partial \mathbb{D}$ to $y \in \partial \mathbb{D}$ if
$\varphi(x) = \varphi(y) \in \mathcal T_n \setminus e_0$. 
 If $\gamma$ does not pass through any edge in $\mathcal T_n^{-}(e_0)$, then the radial projection of $\varphi^{-1}(\gamma)$ onto  $\partial\D$ contains $S^{-}(e_0)$ and the metric
$\rho_0=  {\bf 1}_{A(0;1,2)}$ assigns length $\ge s^{-}(e_0)$ to $\varphi^{-1}(\gamma).$

 In general, the inclusion
$$
S^{-}(e_0) \subset \pi_{\rad}(\varphi^{-1}(\gamma)) \cup \bigcup_{\substack{e \in \mathcal T_n^{-}(e_0) \\ \gamma \cap e \ne \emptyset}} \mathcal \varphi^{-1}(\mathcal T(e))
$$
shows that
$$
\int_{\varphi^{-1}(\gamma)}\Bigl(  {\bf 1}_{A(0;1,2)} + \sum_{e\in\mathcal T_n^-(e_0)} s(e)\rho_e\Bigr)|dz|   \geq s^{-}(e_0),
$$
which proves the admissibility of $\rho$. Together with \eqref{eq:area}, this shows the upper bound
$$
M(\Gamma_{\circlearrowleft}(A)) \, \leq \, A(\rho) \, \lesssim \, \frac{1}{s^-(e_0)^2} \biggl [  1 + \sum_{e \in \mathcal T_n^{-}(e_0)} s(e)^2 \biggr ].
$$
We have thus proved the following theorem:

\begin{theorem}
\label{edges-do-not-shrink-thm}
  Suppose $\{\mathcal T_n\}_{n=0}^\infty$ is an infinite sequence of conformally balanced trees whose vertices have uniformly bounded degrees, for which the sums $S_n = \sum_{e \in \mathcal T_n} s(e)^2$ are uniformly bounded. If  $e_0^{(n)} \subset \mathcal T_n$ is a sequence of edges with $\inf s(e_0^{(n)}) > 0$, then any subsequential Hausdorff limit of the edges $e_0^{(n)}$ is a real-analytic arc.
\end{theorem}

We now apply the above theorem to the sequence of the finite truncations $\{\mathcal T_n\}$ of the infinite trivalent tree.
Inspection shows that
$$
    s(e) \asymp 2^{-d_{\mathcal T_n}(v_{\Root}, e)}.
$$
  As the number of edges $v \in \mathcal T_n$ with $d_{\mathcal T_n}(v_{\Root}, e) = m$ is $\asymp 2^m$, the sums
  $$
    S_n = \sum_{e \in \mathcal T_n} s(e)^2
  $$
  are uniformly bounded in $n = 1, 2, \dots$. If $e \subset \mathcal T$ is an edge in the infinite tree, then its representative $e^{(n)} \subset \mathcal T_n$ has $s(e^{(n)}) \asymp 2^{-d_{\mathcal T_\infty}(v_{\Root}, e) }$. By the theorem above, the diameters of the edges $e^{(n)}$ are bounded from below.

\begin{remark}\label{rmk:diameters} 
If $\mathcal T_n$ is a {\it random} conformally balanced trivalent tree with $n$ edges, chosen uniformly among all of them, then it is not hard to show that the expectation 
$$E[S_n] = \sum_{e \in \mathcal T_n}E[ s(e)^2]  $$
tends to $\infty$ as $n\to\infty,$ suggesting that the diameters of the edges tend to zero. Indeed, it is known \cite{lin-rohde} that the diameters tend to zero with a power of $1/n$, with high probability.
\end{remark}

\section{Structure of a subsequential limit}
\label{sec:limit-structure}

Let $\mathcal T_n$ be the conformally balanced trivalent tree of depth $n$.
In this section, we show that any subsequential limit of the $\mathcal T_n$ has the right topological type:

\begin{theorem}
  \label{top-type}
  For any subsequential Hausdorff limit of the $\mathcal T_n$, one can find a homeomorphism of the plane which takes it onto the Farey tree $\mathcal F$ of Section \ref{sub:Farey}
union the unit circle $\partial \mathbb{D}$.
\end{theorem}

We first pass to a subsequence so that every edge in the infinite trivalent tree has a limit along this sequence. In the previous section, we saw that the limit of each edge is a real-analytic arc. We write $\mathcal T_\infty$ for the union of the Hausdorff limits of the individual edges.
We pass to a further subsequence so that the finite trees $\mathcal T_n$ also possess a Hausdorff limit, which we denote by $\mathcal T_\infty \sqcup \Lambda$. We refer to $\Lambda$ as the {\em limit set}.

The proof of Theorem \ref{top-type} is based on a number of moduli estimates, which control the geometry of the finite trees $\mathcal T_n$.
With the help of these moduli estimates, we prove the following assertions:

\begin{enumerate}[leftmargin=1.75cm, label={\rm (SL\arabic*)}]
\item \label{item:subsequential-limit1} $\mathcal T_\infty$ is dense in the Hausdorff limit of the finite trees $\mathcal T_n$.
\item \label{item:subsequential-limit2} For any branch $[v_0, v_1, v_2, v_3, \dots]$ of $\mathcal T_\infty$ with $d_{\mathcal T_\infty}(v_m, v_{\Root}) = m$, $\lim_{m \to \infty} v_m$ exists.
\item \label{item:subsequential-limit3} Given two branches $[v_0, v_1, v_2, v_3, \dots]$, $[w_0, w_1, w_2, w_3, \dots]$, $\lim_{m \to \infty} v_m = \lim_{m \to \infty} w_m$ if and only if the limits of the corresponding branches in the Farey tree are the same.
\end{enumerate}

We then show the following two topological assertions:
 
 \begin{enumerate}[leftmargin=1.75cm, label={\rm (SL\arabic*)}]
 \setcounter{enumi}{3}
\item \label{item:subsequential-limit4} The limit set $\Lambda$ is a Jordan curve $\partial \Omega$ which encloses $\mathcal T_\infty$.
\item \label{item:subsequential-limit5} There is a natural correspondence between the complementary regions of $\mathcal T_\infty \cup \partial \Omega$ and $\mathcal F \cup \partial \mathbb{D}$.
\end{enumerate}

 From here, the proof of Theorem \ref{top-type} runs as follows: 
 
 \begin{proof}[Proof of Theorem \ref{top-type}]
 Let $h$ be a homeomorphism of $\mathcal T_\infty$ onto the Farey tree $\mathcal F$, which takes vertices to the corresponding vertices. The above properties imply that $h$ extends to a homeomorphism of the closures: $\mathcal T_\infty \cup \partial \Omega$ and $\mathcal F \cup \partial \mathbb{D}$. Since the complementary regions are Jordan domains, we can extend $h$ to a homeomorphism of the plane.
\end{proof}

\subsection{Shrinking of diameters}
\label{sec:shrinking-diameters}

For a vertex $v \in \mathcal T_n$, we denote the subtree which consists of $v$ and its descendants by $\mathcal T_n(v)$.
To prove \ref{item:subsequential-limit1} and \ref{item:subsequential-limit2}, we show:

\begin{lemma}
  \label{diameters-shrink}
  {\em (i)} The diameters of $\mathcal T_n(v)$ tend to zero as $d_{\mathcal T_n}(v_{\Root}, v) \to \infty$, uniformly in  $n$.

    {\em (ii)} The diameters of $\mathcal T_n(vLR^k) \cup  \mathcal T_n(vRL^k)$ tend to zero if either  $d_{\mathcal T_n}(v_{\Root}, v) \to \infty$ or $k \to \infty$, again uniformly in $n$.
\end{lemma}

As in the case of the Farey tree $\mathcal F$ in the unit disk, the diameter of $\mathcal T_n(v)$ depends on the nature of the word representing $v$. If the path joining $v_{\Root}$ to $v$ switches between left and right turns regularly, then the diameters of $\mathcal T_n(v)$ decrease exponentially quickly. On the other hand, if the word for $v$ has long sequences of consecutive $L$'s and $R$'s, then the diameters of $\mathcal T_n(v)$ shrink at a polynomial rate. This dichotomy is reflected in the two types of estimates below.

\paragraph*{Hyperbolic decay.}
\label{sec:hyp-decay}
At an interior vertex $v \in \mathcal T_n$, the domain $\sph \setminus \mathcal T_n$ has three prime ends. Assuming that $v \ne v_{\Root}$ is not the root vertex, we can name the three prime ends as left, right and middle. The {\em left} prime end lies between
$\overline{v_{\parent}v}$ and $\overline{vv_L}$, while the {\em right} prime end lies between $\overline{v_{\parent}v}$ and $\overline{vv_R}$.
Naturally, the {\em middle} prime end lies between $\overline{vv_L}$ and $\overline{vv_R}$. 

\begin{figure}[h]
\centering
    \includegraphics[scale=1]{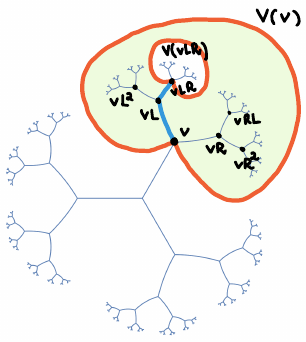}
    \caption{To a non-root vertex $v \in \mathcal T_n$, we associate the domain $V(v)$, bounded by the curve $\gamma(v)$.}
    \label{fig:Vregions}
\end{figure}

Let $\gamma(v)$ denote the hyperbolic geodesic in $\hat{\mathbb{C}} \setminus \mathcal T_n$ which joins the left and right prime ends at $v$ and $V(v)$ be the domain enclosed by $\gamma(v)$, see Fig.~\ref{fig:Vregions}.
With this definition, a vertex $w$ is contained in $V(v)$ if and only if $w$ is represented by a word which begins with $v$.
Moreover, if $v_2$ is a descendant of $v_1,$ then
$V(v_2) \subset V(v_1)$.

\begin{lemma}
  \label{hyp-decay}
  Suppose $v$ is an interior vertex of $\mathcal T_n$, other than the root vertex. Then, $\Mod V(v) \setminus V(vLR) \asymp 1$ and $\Mod V(v) \setminus V(vRL) \asymp 1$.
\end{lemma}

It is enough to show the statement regarding $\Mod \Ann$ for $\Ann = V(v) \setminus V(vLR)$ as the situation with $\Mod V(v) \setminus V(vRL)$ is entirely symmetric.
To prove the lemma, we need to give uniform upper bounds for the moduli of the curve families $\Gamma_{\circlearrowleft}(\Ann)$ and $\Gamma_{\uparrow}(\Ann)$,
which are independent of $n$ and $v \in \mathcal T_n$. 

To deal with the first curve family, simply note that every $\gamma\in\Gamma_{\circlearrowleft}(\Ann)$ intersects at least one of the two edges $\overline{v v_L}$ and
$\overline{v_L v_{LR}}$ so that the sum of the two obstacles $\rho=\rho_{\overline{v v_L}} + \rho_{\overline{v_L v_{LR}}}$ is an admissible metric of area $A(\rho)=O(1).$

To deal with the second curve family, by conformal invariance, we may give an upper bound for the modulus of the curve family $\varphi^{-1}(\Gamma_\uparrow(\Ann))$ 
in $\varphi^{-1}(\Ann) \subset \sph \setminus \mathbb{D}$ which allows teleportation, as we did before in Section \ref{sec:interior-edges}.
Cutting $\Ann$ along the tree, we obtain a conformal rectangle $\Rect = \Ann \setminus \mathcal T_n$  whose vertices are the prime ends where $\gamma(v)$ and $\gamma(vLR)$ meet $\mathcal T_n$. Its pre-image $\hat\Rect = \varphi^{-1}(\Rect) \subset \sph \setminus \mathbb{D}$ is a conformal rectangle
whose vertices are the points where the geodesics $\varphi^{-1}(\gamma(v))$  and  $\varphi^{-1}(\gamma(vLR))$ meet the unit circle. We label the vertices $z_1, z_2, z_3, z_4$ in
counter-clockwise order such that $z_1$ corresponds to the right prime end of $v$. Due to the ``left-right'' turn between $v$ and $vLR$, the distances between the points $z_i$, $1 \le i \le 4$, are comparable:
$$
|z_1 - z_2| \asymp |z_2 - z_3| \asymp |z_3 - z_4| \asymp 2^{-d}, \qquad d = d_{\mathcal T_n}(v_{\Root},v).
$$
Hence, the background metric  $\rho_0={\bf1}_{\varphi^{-1}(\Ann)}$ assigns length $\gtrsim 2^{-d}$ to every curve in $\varphi^{-1}(\Gamma_\uparrow(\Ann))$ that does not teleport.

Arguing as in Section \ref{sec:interior-edges} shows that the metric
\begin{equation}
\label{eq:testmetric}
\rho = C_0 2^d \Bigl( \rho_0 + \sum_{e\in V(v)} s(e)  \rho_e \Bigr)
\end{equation}
is admissible if $C_0$ is sufficiently large (independent of $n$ and $v$). 
More precisely, while the set  $\varphi^{-1}(\gamma)$ may be disconnected,
$$
\sigma=\varphi^{-1}(\gamma)\cup\bigcup_{e\cap\gamma\neq\emptyset} S(e)
$$
is connected and intersects both geodesics $\varphi^{-1}(\gamma(v))$  and  $\varphi^{-1}(\gamma(vLR))$. Inspection shows that the integral $\int_\sigma \rho_0 |dz|$ computes the Euclidean length of $\sigma\setminus\partial\D$, whereas $\int_\sigma \bigl ( \sum_{e\in V(v)} s(e)  \rho_e \bigr ) |dz|$ is bounded below by the Euclidean length of $\sigma\cap\partial\D$. 
As a result,
$$
\int_{\sigma} \biggl \{ \rho_0 + \sum_{e\in V(v)} s(e)  \rho_e \biggr \} |dz|
$$ is greater or equal to the Euclidean distance between the geodesics  $\varphi^{-1}(\gamma(v))$  and  $\varphi^{-1}(\gamma(vLR))$, which is comparable to $2^{-d}$. Consequently, the factor $C_0 2^d$ in (\ref{eq:testmetric}) makes the metric $\rho$ admissible.

From $s(e) \asymp 2^{-d_{\mathcal T_n}(v_{\Root}, e)}$, it is clear that $\sum_{e\in V(v)} s(e)^2 \lesssim 2^{-2d}$. The area bound $A(\rho)=O(1)$ now follows from \eqref{eq:area}. Putting the above information together shows the desired modulus bound.

\begin{figure}[h]
\centering
    \includegraphics[scale=1]{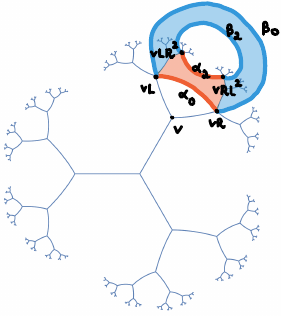}
    \caption{The region $W_j = W_j(v)$ is associated to a non-root vertex $v \in \mathcal T_n$ and an integer $j \ge 0$. It is bounded by the hyperbolic geodesics $\alpha_j$ and $\beta_j$.}
\label{fig:Wregions}
\end{figure}

\paragraph*{Parabolic decay.}
\label{sec:par-decay}
We continue to assume that $v \in \mathcal T_n$ is an interior vertex, other than the root vertex. For each $0 \le j \le n - 1 - d_{\T_n}(v_{\Root}, v)$, we connect the vertices $vLR^j$ and  $vRL^j$ by two hyperbolic geodesics $\alpha_j, \beta_j \subset \sph \setminus \mathcal T_n$, with the {\em inner geodesic}\/ $\alpha_j$ joining
$$
(vLR^j)_{\Right} \quad \text{with} \quad (vRL^j)_{\Left}
$$
and the {\em outer geodesic}\/  $\beta_j$ joining
$$
(vLR^j)_{\Left} \quad \text{with} \quad  (vRL^j)_{\Right}.
$$
We then define $W_j = W_j(v)$ as the simply-connected domain bounded by the Jordan curve $\alpha_j \cup \beta_j$.
See Fig.~\ref{fig:Wregions}.

\begin{lemma}
  \label{par-decay}
  Suppose $v$ is an interior vertex of $\mathcal T_n$, other than the root vertex. Then, $$\Mod \, W_0(v) \setminus W_k(v) \asymp \log(1+k),$$ for any $1 \le k \le n - 1 - d_{\T_n}(v_{\Root}, v)$.
\end{lemma}

Since the annulus $V(v) \setminus V(vLR^k) \supset W_0 \setminus W_k$, its modulus is strictly larger. In particular, the lemma implies
that $\Mod V(v) \setminus V(vLR^k) \gtrsim \log(1+k)$.

\begin{proof}
  For brevity, we write $\Ann = W_0 \setminus W_k$. To show the  upper bound for $\Mod \Ann$, we need to estimate the modulus of the family of curves $\Gamma_{\circlearrowleft}$ which separate the boundary components of $\Ann$.
  The tree $\mathcal T_n$ splits $\Ann$ into two conformal rectangles $\Rect_\alpha$ and $\Rect_\beta$, with  $\pRect_\alpha \supset \alpha_0 \cup \alpha_k$ and $\pRect_\beta \supset \beta_0 \cup \beta_k$. Since a curve in  $\Gamma_{\circlearrowleft}(\Ann)$ contains a crossing that joins the $\mathcal T_n$-sides of $\Rect_\alpha$, $\Mod \Gamma_{\circlearrowleft}(\Ann) \le \Mod \Rect_\alpha$. The latter modulus may be computed in the exterior unit disk: $\Mod \Rect_\alpha = \Mod \varphi^{-1}(\Rect_\alpha) \asymp \log(1+k)$ as desired.

  We now turn to the lower bound. For this purpose, we decompose $\Ann$ into a union of shells:
  $$
    \Ann \, = \, \bigcup_{j=1}^{k} \Ann_j \, = \, \bigcup_{j=1}^{k} W_{j-1} \setminus W_{j}.
  $$
By the parallel rule, it is enough to show that $\Mod \Gamma_{\uparrow}(\Ann_j) \lesssim j$, for each $j = 1, 2, \dots, k$. As usual, we estimate the modulus of the  family $$\varphi^{-1}(\Gamma_{\uparrow}(\Ann_j)) \subset \varphi^{-1}(\Ann_j) \subset \sph \setminus \mathbb{D}.$$
  The pre-image $\varphi^{-1}(\Ann_j) = \hat\Rect_{\beta, j} \cup \hat\Rect_{\alpha, j}$ consists of two conformal rectangles in $\sph \setminus \mathbb{D}$, with $\hat\Rect_{\alpha, j}$ bounded by $\hat{\alpha}_{j-1}, \hat{\alpha}_{j}$ and the unit circle, and
  $\hat\Rect_{\beta, j}$ bounded by $\hat{\beta}_{j-1}, \hat{\beta}_{j}$ and the unit circle.

Let $\rho_{\alpha, 0}(z)$ be the extremal metric on the conformal rectangle $\hat\Rect_{\alpha, j}$ for the family of curves contained in  $\hat\Rect_{\alpha, j}$ that connect $\hat{\alpha}_{j-1}$ and $\hat{\alpha}_{j}$.
It is easy to see that $A(\rho_{\alpha, 0}) \asymp j+1.$
As in the proof of Lemma \ref{hyp-decay}, there is a metric $\rho_{\beta, 0}$ of the form \eqref{eq:testmetric} with $A(\rho_{\beta, 0})\asymp 1$ which  assigns length $\geq 1$ to every curve $\gamma$ in $\hat\Rect_{\beta, j}$ that connects $\hat{\beta}_{j-1}$ and $\hat{\beta}_{j}$ with or without teleportation. More precisely, since the four marked endpoints of $\hat{\beta}_{j-1}$ and $\hat{\beta}_{j}$ have mutually comparable distances $\asymp 2^{- d_{\T_n}(v_{\Root}, v) - j}$, the reasoning in the proof of 
Lemma \ref{hyp-decay} shows that the metric 
\begin{equation}
\label{eq:testmetric-par}
\rho_{\beta, 0} = C_0 2^{d_{\T_n}(v_{\Root}, v) + j} \biggl( \rho_0 + \sum_{e\in V(v_{LR^{j-1}})\cup V(v_{RL^{j-1}})} s(e)  \rho_e \biggr)
\end{equation}
is admissible where $\rho_0={\bf1}_{\hat\Rect_{\beta, j}}$ and $C_0$ is sufficiently large.

A path in $\varphi^{-1}(\Gamma_{\uparrow}(\Ann_j))$ connects $ \hat{\alpha}_{j-1} \cup \hat{\beta}_{j-1}$ with $\hat{\alpha}_{j} \cup \hat{\beta}_{j}$, where one is allowed to take shortcuts by
teleporting from $x \in  \partial \mathbb{D}$ to $y \in \partial \mathbb{D}$ if
$\varphi(x) = \varphi(y) \in \mathcal T_n$. Such a path is
either contained in $\hat\Rect_{\alpha, j}$, or contained in $\hat\Rect_{\beta, j}$, or intersects one of the two edges  $e_1=[vLR^{j-1}, vLR^{j}]$ and $e_2=[vRL^{j-1}, vRL^{j}].$ 
To obtain a metric admissible for $\varphi^{-1}(\Gamma_{\uparrow}(\Ann_j))$, we modify $\rho_0(z) = \rho_{\alpha, 0}(z) + \rho_{\beta, 0}(z)$ by adding two obstacles along the edges $e_1$ and $e_2$ which make it impractical for a path to teleport from $\hat\Rect_{\beta, j}$ to $\hat\Rect_{\alpha, j}$:
  $$
    \rho = (\rho_{\alpha, 0} + \rho_{\beta, 0}) + \rho_{e_1} + \rho_{e_2}.
  $$
As each obstacle has area $O(1)$, the area $A(\rho) \asymp j+1$, which gives the desired modulus bound.
\end{proof}

\paragraph*{Putting this together.} We are now ready to show Lemma \ref{diameters-shrink}:
  
\begin{proof}[Proof of Lemma \ref{diameters-shrink}]
Let $v \in \mathcal T_n$ be an interior vertex, other than the root vertex. From the definitions, it is clear that
$
  \mathcal T_n(v) \subset V(v).
$
  Let $$[v_{\Root}, v] = [v_0 = v_{\Root}, v_1, v_2, v_3, \dots, v_m = v]$$ be the path in $\mathcal T_n$ joining $v_{\Root}$ to $v$.
  In view of the hydrodynamic normalization, $V(v_1) \subset \Omega_2 \subset B(0,8)$ is contained in a ball of fixed size. Consequently, to prove (i), it is enough to show that
  $\Mod V(v_1) \setminus V(v)$ is large when $d_{\mathcal T_n}(v_{\Root}, v)$ is large.

  There are two cases to consider. If the path $[v_{\Root}, v]$ frequently switch between left and right turns, then
  $\Mod V(v_1) \setminus V(v)$ will be large by Lemma \ref{hyp-decay} and the parallel rule. If we turn left many times or turn right many times without switching, then  $\Mod V(v_1) \setminus V(v)$ will be large by Lemma \ref{par-decay}. In both cases, $\diam V(v) \to 0$ uniformly in $n$ as $d_{\mathcal T_n}(v_{\Root}, v)\to\infty$.
To prove (ii), we note that 
$$\mathcal T_n(vLR^k) \cup \mathcal T_n(vLR^k) \, \subset \, W_k(v) \, \subset \, V(v)$$ 
and appeal to Lemma \ref{par-decay}.
\end{proof}

\subsection{The limit set is a Jordan curve}

Our next objective is to show \ref{item:subsequential-limit3} and \ref{item:subsequential-limit4}. For two vertices $v_1, v_2 \in \mathcal T_n$, we denote by $d_\omega(v_1, v_2)$ the harmonic measure as seen from infinity of the shortest arc on the unit circle that contains a point of $\varphi^{-1}(v_1)$ and a point of $\varphi^{-1}(v_2)$.  The following lemma says that if the harmonic measure between two boundary vertices $v_1,v_2$ is small, then the Euclidean distance $|v_1-v_2|$ is also small:

\begin{lemma}
  \label{does-not-get-torn-apart}
  For any $\varepsilon > 0$, there exists an $\eta > 0$, such that if  $v_1, v_2 \in \mathcal T_n$ are two boundary vertices for which $d_\omega(v_1, v_2) < \eta$, then
  the Euclidean distance $|v_1-v_2| < \varepsilon$.
\end{lemma}

We explain the proof via an analogy: if $x_1, x_2$ are two points on the unit circle, then either $x_1, x_2$ are contained in a single  dyadic arc  of length comparable to $|x_1 - x_2|$ or they are contained in the union of two adjacent dyadic arcs whose lengths  are comparable to $|x_1 - x_2|$.
Similarly, in the trivalent tree, one has two non-mutually exclusive possibilities: Denote $v$ the last common ancestor of $v_1$ and $v_2$ so that $v_1, v_2 \in \mathcal T_n(v)$ and $v_1=vLX, v_2=vRY$ (or vice versa) for some sequences $X,Y$. Then at least one of the following statements is true:
\begin{enumerate}
\item $\omega_{\sph \setminus \mathcal T_n, \infty}(\mathcal T_n(v)) \asymp d_{\omega}(v_1,v_2)$,
\item For the maximal integer $k \ge 1$ so that $v_1 \in \mathcal T_n(vLR^k)$ and $v_2 \in \mathcal T_n(vRL^k)$, we have
$\omega_{\sph \setminus \mathcal T_n, \infty} \bigl (\mathcal T_n(vLR^k) \cup  \mathcal T_n(vRL^k) \bigr ) \asymp d_{\omega}(v_1,v_2)$.
\end{enumerate}
In either case, one may use Lemma \ref{diameters-shrink} to show that the Euclidean distance $|v_1 - v_2|$ is small if  $d_\omega(v_1, v_2)$ is small.

The same argument shows that for any $\varepsilon > 0$, there exists an $\eta > 0$, such that for any two boundary vertices  $v_1, v_2 \in \mathcal T_n$ with $d_\omega(v_1, v_2) < \eta$, the union of the geodesics that join consecutive boundary vertices of $\T_n$ between $v_1$ and $v_2$ has diameter $< \varepsilon$. Indeed, these geodesics are contained in the region $V(v)$ in the first case above and in the region $W_k(v)$ in the second case.

We now show the converse to Lemma \ref{does-not-get-torn-apart}, namely, if the harmonic measure between two boundary vertices in  $\mathcal T_n$ is bounded below, then so is their Euclidean distance:

\begin{lemma}
  \label{does-not-collapse}
  For any $\varepsilon > 0$, there exists an $\eta > 0$, such that if  $v, w \in \mathcal T_n$ are two boundary vertices for which $d_\omega(v, w) > \eta$, then
  the Euclidean distance $|v-w| > \varepsilon$.
\end{lemma}

\begin{proof}
  Let $[v_0 = v_{\Root}, v_1, v_2, v_3, \dots, v_n = v]$ be the path in $\mathcal T_n$ joining $v_{\Root}$ to $v$ and
  $[w_0 = v_{\Root}, w_1, w_2, w_3, \dots, w_n = w]$ be the path joining $v_{\Root}$ to $w$.
  The assumption implies that there exists an $n_0 = n_0(\eta) \ge 1$ sufficiently large so that
  the harmonic measure between $E = [v_{n_0},v]$ and $F = [w_{n_0}, w]$ is at least $\eta/2$.

  Recall that in Section \ref{sec:interior-edges}, we showed that the diameters of  $E$ and $F$ are bounded from below. To show that $E$ and $F$ are a definite distance apart, it is enough to give an upper bound for the modulus of the family of curves $\Gamma_{E \leftrightarrow F}$ that connect $E$ to $F$ in $\Omega_2$. By the conformal invariance, we may instead estimate the modulus $\varphi^{-1}(\Gamma_{E \leftrightarrow F})$
  in $A(0;1,2)$ where teleportation is allowed between the pre-images of points in $\mathcal T_n$.

An argument similar to the one in Section \ref{sec:interior-edges} shows that
$$\int_{\varphi^{-1}(\gamma)}\Bigl(  {\bf 1}_{A(0;1,2)} + \sum_{e\in\T_n} s(e)\rho_e\Bigr)|dz|   \geq \eta/2, \qquad \gamma \in \Gamma_{E \leftrightarrow F},
$$
i.e.~$2/\eta$ times the integrand is an admissible metric $\rho$ with $A(\rho)= O(1/\eta^2).$
\end{proof}

Lemmas \ref{does-not-get-torn-apart} and \ref{does-not-collapse} imply that $\Lambda = \partial \Omega$ is a Jordan curve: Indeed, joining consecutive boundary vertices of $\T_n$ by hyperbolic geodesics, we obtain a sequence of Jordan curves $\Lambda_n.$ If we parametrize these curves by harmonic measure from infinity, then they converge uniformly to a continuous limit curve by Lemma \ref{does-not-get-torn-apart}, and this limit curve is simple by Lemma \ref{does-not-collapse}. Furthermore, it is disjoint from $\mathcal T_\infty$ by Lemma \ref{basic-trees1}.

\subsection{Formation of \texorpdfstring{$\Omega$-horoballs}{Ω-horoballs}}

We now turn to showing \ref{item:subsequential-limit5}. Let $v_0\ne v_{\Root}$ be a vertex of the infinite trivalent tree.  
For $j \ge 1$, set
$$
  v_j = v_0 LR^{j-1} \qquad \text{and}
  \qquad v_{-j} = v_0 RL^{j-1}.
$$
From  \ref{item:subsequential-limit3}, the limits
$$
\lim_{j \to +\infty}  v_j \quad \text{and} \quad \lim_{j \to -\infty}  v_{j}
$$
exist and are equal.
We refer to their common value $p$ as a {\em cusp} or {\em parabolic point}\/. In particular, the union of the edges
$$
  \bigcup_{j=-\infty}^\infty \overline{v_jv_{j+1}} \, \subset \, \mathcal T_\infty
$$
defines a Jordan curve, which we denote $\partial \Omega_p$. At the root vertex, one can similarly construct three Jordan domains $\Omega_{p_1}, \Omega_{p_2}, \Omega_{p_3}$.
We refer to the regions  $\{\Omega_{p_i}\}$ as $\Omega$-horoballs.

\begin{lemma}
  \label{bounded-components}
  The regions $\{\Omega_{p_i}\}$ enumerate the bounded components of $\mathbb{C} \setminus \lim \mathcal T_n$. 
\end{lemma}

\begin{proof}
  We approximate the regions $\Omega_{p_i}$ by Jordan domains $\Omega_{p_i}^{(n)}$ constructed using the finite approximating trees $\mathcal T_n$ as follows:
  Each finite tree $\mathcal T_n$ contains only finitely many corresponding vertices $\{v_j\}_{j=-m}^m$, where $m = n - d(v_{\Root}, v_0)$. The union of the edges $\bigcup_{j=-m}^{m-1} \overline{v_jv_{j+1}} \subset \mathcal T_n$ is a Jordan arc. To form $\partial \Omega_{p_i}^{(n)}$, we close this Jordan arc with the hyperbolic geodesic $\alpha_{p_i}^{(n)} \subset \sph \setminus \mathcal T_n$ that connects the leaves $v_{-m}, v_m \in \T_n$.
  
  In view of Lemma \ref{par-decay}, $\diam \alpha_{p_i}^{(n)} \to 0$ and  $\Omega_{p_i} = \lim \Omega_{p_i}^{(n)}$. Since $\Omega_{p_i}^{(n)}$ is disjoint from the tree $\mathcal T_n$, the regions $\Omega_{p_i} = \lim \Omega_{p_i}^{(n)}$ are indeed bounded  components of the complement $\mathbb{C} \setminus \lim \mathcal T_n$.

 Can there be any more complementary components? If $O$ is any connected component of $\mathbb{C} \setminus \lim \mathcal T_n$, then $\partial O\subset \T_{\infty}\cup \Lambda.$  If $\partial O$ intersects one of the edges of $\T_{\infty},$ then $O$ is one of the four $\Omega$-horoballs who form a neighborhood of this edge. If $\partial O$ does not intersect $\T_{\infty},$ then 
$\partial O\subset\Lambda$, and since $\Lambda$ is a Jordan curve, $O$ must be the unbounded component of $\mathbb{C} \setminus \Lambda$.
\end{proof}

Having established Properties \ref{item:subsequential-limit1}--\ref{item:subsequential-limit5}, the proof of Theorem \ref{top-type} is complete.

\subsection{Uniqueness of the limit}
\label{sec:uniqueness-estimate}

For a non-root vertex $v \in \mathcal T$, we define the {\em shadow} $s_v \subset \partial \Omega$ as the shorter arc of $\partial \Omega$ which joins
$vLRL^\infty = \lim_{m \to \infty} vLRL^m$ and $vRLR^\infty = \lim_{m \to \infty} vRLR^m$.
A brief inspection of the homeomorphic picture of the Farey tree $\mathcal F \subset \mathbb{D}$ shows that any point on $\partial \Omega$ that is not a cusp of an $\Omega$-horoball is contained in infinitely many shadows.
The following estimate will be used in Section \ref{sec:uniqueness} in conjunction with Lemma \ref{weak-conformal-removability} to show that the Hausdorff limit of the true trees $\mathcal T_n$ is unique:

\begin{lemma}
\label{tree-square-sum}
The sums
  \begin{equation}
    \label{eq:tree-square-sum}
    \sum_{v \in \mathcal  T_n, \, v \ne v_{\Root}}  \Bigl \{ \diam^2 V(vRL) + \diam^2 V(vLR) \Bigr \}
  \end{equation}
  are uniformly bounded above, independent of $n$.
\end{lemma}

Since $s_v$ is the Hausdorff limit as $n \to \infty$ of $(V(vRL) \cup V(vLR)) \cap \partial \Omega$, the above lemma implies that $ \sum_{v \in \mathcal  T_\infty, \, v \ne v_{\Root}} \diam^2 s_v < \infty$. In particular, $\partial \Omega$ has area zero.

\begin{proof}
  For a hyperbolic geodesic $\hat \gamma \subset \{z \in \mathbb{C} : 1 < |z| < 2\} \subset \hat{\mathbb{C}} \setminus \mathbb{D}$, let $z_{\hat \gamma}$ be the Euclidean midpoint of $\hat \gamma$
  and $B_{\hat \gamma} \subset \hat{\mathbb{C}} \setminus \mathbb{D}$ be the ball of hyperbolic radius 1/10 centered at $\frac{1 + |z_{\hat \gamma}|}{2} \cdot z_{\hat \gamma}$. In view of the restriction on $\hat\gamma$, the ball $B_{\hat \gamma}$ is contained in the bounded domain enclosed by  $\hat \gamma$ and the unit circle.

  Similarly, to a hyperbolic geodesic $\gamma \subset \Omega_2 \subset \hat{\mathbb{C}} \setminus \mathcal T_n$, we can associate the topological disk $B_\gamma := \varphi(B_{\varphi^{-1} \gamma})$.
  By Koebe's distortion theorem, $B_\gamma$ is approximately round in the sense that its area is comparable to its diameter squared.
  
We apply the above construction to the geodesics $\gamma(v)=\partial V(v) \subset \Omega_2$ from Section \ref{sec:shrinking-diameters}, where $v \ne v_{\Root}$ ranges over interior vertices of $\mathcal T_n$, other than the root vertex. From the construction, it is clear that $B_{\gamma(v)} \subset V(v)$.

To prove the estimate (\ref{eq:tree-square-sum}), it is enough to show that
  \begin{equation}
    \label{eq:uniqueness-estimate}
    \diam  V(vLR)
    \asymp \diam B_{\gamma(vLR)},
  \end{equation}
  as the topological disks $B_{\gamma(vLR)}$ are disjoint and are contained in a bounded set. In view of Lemma \ref{lem:egg_yolk}, to prove (\ref{eq:uniqueness-estimate}), we may show the following two moduli estimates:
  \begin{enumerate}
    \item  $\Mod V(v) \setminus V(vLR)$ is bounded below.
    \item $\Mod V(v) \setminus B_{\gamma(vLR)}$ is bounded above.
  \end{enumerate}
  The first estimate was already established in Lemma \ref{hyp-decay}. The second estimate is automatic from Koebe's distortion theorem.
\end{proof}

\begin{remark}
Let $\Omega \subset \mathbb{C}$ be a Jordan domain, $K \subset \Omega$ be a compact set and  $z_0 \in \Omega$ be an interior point.
In the work of Jones and Smirnov \cite{jones-smirnov}, the {\em shadow} of $K$ with respect to $z_0 \in \Omega$ is defined as the set of endpoints of hyperbolic rays emanating from $z_0$ which pass through $K$. It is not difficult to show that the set $s_v$ described above and the Jones-Smirnov shadow of the closed ball of hyperbolic radius 1 centered at $v$ with respect to $v_{\Root} \in \Omega$ intersect and have comparable diameters.
\end{remark}

\section{Convergence}
\label{sec:convergence}

In this section, we show that the Hausdorff limit $\mathcal T_\infty \cup \partial \Omega$ of the finite trees $\mathcal T_n$ is unique. The main step is to prove that it realizes the mating of $z \to \overline{z}^2$, acting on the exterior unit disk $\mathbb{D}_e$, and the Markov map $\rho(z)$ associated to the reflection group of an ideal triangle $\triangle_{\hyp}$, acting on the unit disk $\mathbb{D}$.

\subsection{Relative harmonic measure}
\label{sec:relative}

Suppose that $U$ is a Jordan domain and  $p \in \partial U$.  While it does not make sense to talk about the harmonic measure of an arc $I \subset \partial U$ as viewed from $p$,
one can talk about the {\em relative harmonic measure} of two arcs $I, J \subset \partial U$ that do not contain $p$\,:
$$
  \omega_{U, p}(I, J) = \lim_{z \to p} \frac{ \omega_{U, z} (I)}{ \omega_{U, z} (J)}.
$$
It is easy to see that the quantity $\omega_{U, p}(I, J)$ varies continuously provided that $p$ stays away from $I \cup J$. More precisely, if a sequence of Jordan quadruples $(U_n, p_n, I_n, J_n)$ converges in the Hausdorff topology to a Jordan quadruple $(U, p, I, J)$, then
$$
  \omega_{U_n, p_n}(I_n, J_n) = \lim_{n \to \infty}\omega_{U, p}(I, J).
$$

\begin{example}
  When $\Omega = \mathbb{H}$, $p = \infty$ and $J=[0,1]$, the relative harmonic measure $\omega_{\mathbb{H}, \infty}(\cdot, [0,1])$ is just Lebesgue measure on the real line.
\end{example}

\subsection{Farey horoballs}
\label{sec:farey-horoballs}

The Farey tree $\mathcal F$ partitions the unit disk $\mathbb{D}$ into regions which we call {\em Farey horoballs} $H_{p_i}$. We index the Farey horoballs
by the point where they touch the unit circle. We label the vertices on $\partial H_{p_i}$ in counter-clockwise order by $v^j(H_{p_i})$, $j \in \mathbb{Z}$, with  $v^0(H_{p_i})$ being the vertex with the smallest combinatorial distance to $v_{\Root}$.

By construction, the Farey tree is invariant under the group generated by reflections in the sides of $\triangle_{\hyp}$. As such, Farey horoballs enjoy the
following two properties:

\begin{enumerate}[leftmargin=1.75cm, label={\rm (F\arabic*)}]
  \item \label{item:farey-tree1} Any two edges $e_1, e_2 \subset \partial H_{p_i}$ have the same relative harmonic measure as viewed from $p_i$, i.e.~
        $$
          \omega_{H_{p_i}, p_i}(e_1, e_2)= 1.
        $$
  \item \label{item:farey-tree2}  If an edge $e$ belongs to two neighbouring Farey horoballs $H_{p_i}$ and $H_{p_j}$, then the relative harmonic measures are the same from both sides:
        $$
          \omega_{H_{p_i}, p_i}(I, e)= \omega_{H_{p_j}, p_j}(I, e), \qquad I \subseteq e.
        $$
\end{enumerate}

\subsection{Interior Structure of \texorpdfstring{$\Omega$}{Ω}}
\label{sub:interior}

In Section \ref{sec:limit-structure}, we saw that any Hausdorff limit $\mathcal T_\infty \cup \partial \Omega$ of the $\mathcal T_n$ is ambiently homeomorphic to the Farey tree $\mathcal F$ union the unit circle $\partial \mathbb{D}$.
Recall that the connected components of $\Omega \setminus \mathcal T_\infty$ are called $\Omega$-horoballs and are labeled by the point where they meet $\partial \Omega$.

\begin{lemma}
  The $\Omega$-horoballs also enjoy properties \ref{item:farey-tree1} and \ref{item:farey-tree2}.
\end{lemma}

\begin{proof}
  Since the arguments are very similar, we only present the proof of the second property and leave the proof of the first property to the reader.

  We approximate $\Omega_{p_i}$ by Jordan domains $\Omega_{p_i}^{(n)}$ as in the proof of Lemma \ref{bounded-components}. Pick an arbitrary point $p_i^{(n)} \in \alpha_{p_i}^{(n)}$. As the diameters of $\alpha_{p_i}^{(n)}$ tend to 0, the points $p_i^{(n)} \to p_i$.

  Suppose two neighbouring $\Omega$-horoballs  $\Omega_{p_i}$ and  $\Omega_{p_j}$ meet along an edge $e$. Given an arc  $I \subset e$, we can approximate it in the Hausdorff topology by arcs $I_n \subset e^{(n)} \subset \mathcal T_n$. By the aforementioned continuity of the relative harmonic measure, we have
$\omega_{\Omega_{p_i}, p_i}(I, e)=\lim  \omega_{\Omega_{p_i}^{(n)}, p_{i}^{(n)}}(I_n, e^{(n)})$ so that it is enough to show
$ \omega_{\Omega_{p_i}^{(n)}, p_{i}^{(n)}}(I_n, e^{(n)})\sim  \omega_{\Omega_{p_j}^{(n)}, p_{j}^{(n)}}(I_n, e^{(n)}).$

An intuitive albeit somewhat informal proof of  \ref{item:farey-tree2} is as follows:
  Run Brownian motion from $\infty$ until it hits  $\mathcal T_n$.
  If it is to hit the arc $I_n \subset e^{(n)}$ from the side of $\Omega_{p_i}$, denoted by $I_n \,|\, \Omega_{p_i}^{(n)}$, then it must pass through the gate $\alpha_{p_i}^{(n)}$.
  Since the diameter of the gate $\alpha_{p_i}^{(n)}$ is very small,
  $$
    \omega_{\Omega_{p_i}^{(n)}, p_{i}^{(n)}}(I_n, e^{(n)})
    \, \sim \,
    \frac{ \omega_{ \hat{\mathbb{C}} \setminus  \mathcal T_n, \infty}(I_n \,|\, \Omega_{p_i}^{(n)}) }
    { \omega_{ \hat{\mathbb{C}} \setminus  \mathcal T_n, \infty}(e^{(n)}  \,|\, \Omega_{p_i}^{(n)}) }
    \, = \,
    \frac{ \omega_{ \hat{\mathbb{C}} \setminus  \mathcal T_n, \infty}(I_n  \,|\,  \Omega_{p_j}^{(n)}) }
    { \omega_{ \hat{\mathbb{C}} \setminus  \mathcal T_n, \infty}(e^{(n)}  \,|\,  \Omega_{p_j}^{(n)}) }
    \, \sim \,
    \omega_{\Omega_{p_j}^{(n)}, p_{j}^{(n)}}(I_n, e^{(n)}).
  $$
Here, we have used that the harmonic measures on the two sides of every edge $e$ in a true tree are identical.
  The lemma follows after taking $n \to \infty$.

For a rigorous justification of these asymptotic equalities, notice that the pre-images of the approximate $\Omega$-horoballs $\Omega_{p_i}^{(n)}$ and $\Omega_{p_j}^{(n)}$
under the hydrodynamically normalized Riemann maps
$\varphi_n: \sph \setminus \mathbb{D} \to \sph \setminus \mathcal T_n$ are of the form
$$\varphi_n^{-1}(\Omega_{p_i}^{(n)})= \D_e \cap B_n,$$
where the $B_n$ are (small) disks with centers near $\partial\D_e$ and  $\varphi_n^{-1}(p_i^{(n)}) \in\partial B_n.$
The conformal maps $f_n$ of $B_n\cap\D_e$ onto the upper half-plane $\mathbb{H}$ that send $p_i^{(n)}$ to $\infty$ and $\varphi_n^{-1}(e^{(n)})$ to $[0,1]$ extend by reflection to $B_n$. As the modulus of the annulus $B_n\setminus\varphi_n^{-1}(e^{(n)})$ tends to infinity, by the Koebe distortion theorem,
$$ \omega_{\Omega_{p_i}^{(n)}, p_{i}^{(n)}}(I_n, e^{(n)})=\text{length}(f_n(\varphi_n^{-1}(I_n)))
\sim 
\frac{\text{length}(\varphi_n^{-1}(I_n))}{\text{length}(\varphi_n^{-1}(e^{(n)}))} = 
\frac{ \omega_{ \hat{\mathbb{C}} \setminus  \mathcal T_n, \infty}(I_n \,|\, \Omega_{p_i}^{(n)}) } 
{ \omega_{ \hat{\mathbb{C}} \setminus  \mathcal T_n, \infty}(e^{(n)}  \,|\, \Omega_{p_i}^{(n)}) },
$$
which is what we wanted to show.
\end{proof}

Since $\partial \mathbb{D} \cup \mathcal F$ and $\partial \Omega \cup \mathcal T$ are ambiently homeomorphic, one has a correspondence between the bounded complementary components of
$\partial \mathbb{D} \cup \mathcal F$ (Farey horoballs) and those of
$\partial \Omega\, \cup \mathcal T$ ($\Omega$-horoballs). For each pair of corresponding complementary regions, form the conformal mapping $\varphi_i: H_i \to \Omega_i$ which takes
$$
  p(H_i) \to p(\Omega_i), \quad v^0(H_i) \to v^0(\Omega_i), \quad v^1(H_i) \to v^1(\Omega_i).
$$
(As  Farey horoballs and $\Omega$-horoballs are Jordan domains, by Carath\'eodory's theorem, the maps $\varphi_i$ extend to homeomorphisms between the closures.)

Since Farey and $\Omega$-horoballs possess the property \ref{item:farey-tree1}, $\varphi_i$ maps $v^j(H_i) \to v^j(\Omega_i)$ for any $j \in \mathbb{Z}$.
Moreover, as Farey and $\Omega$-horoballs possess the property  \ref{item:farey-tree2}, we have:

\begin{lemma}[Interior structure lemma]
  \label{inner-structure}
  The mappings  $\varphi_i: H_i \to \Omega_i$  glue up to form a conformal mapping $\varphi: \mathbb{D} \to \Omega.$ In other words, if $H_i$ and $H_j$ share a common edge $e$, then $\varphi_i|_e = \varphi_j|_e$.
\end{lemma}

Indeed, since the edges of $\mathcal F$ are analytic arcs and the homeomorphism $\varphi$ is conformal on $\mathbb{D} \setminus \mathcal F$, it follows that $\varphi$ extends analytically across the open edges. As the vertices are isolated points, they are removable singularities.

\subsection{Exterior Structure of \texorpdfstring{$\Omega$}{Ω}}
\label{sec:exterior-structure}

By definition, the harmonic measure $\omega_{\hat{\mathbb{C}} \setminus \mathcal T_n, \infty}$ is supported on $\mathcal T_n$. From Koebe's $1/4$ theorem, we know that the true trees
$\mathcal T_n \subset B(0, 8)$ are contained in a fixed compact set, so that any subsequential
weak-$*$ limit $\omega$ of the $\omega_{\hat{\mathbb{C}} \setminus \mathcal T_n, \infty}$ is a probability measure supported on  the Hausdorff limit $\mathcal T_\infty \cup \partial \Omega$. As the harmonic measure of any individual edge tends to zero, the support of the limiting measure $\omega$ is contained in $\partial \Omega$.
Finally, since $\partial \Omega$ is uniformly perfect, being a Jordan curve,  $\omega =  \omega_{\Omega_e, \infty}$.

Consider the map $f(z) = \overline{z}^2$ acting on the unit circle. It has fixed points at $1$, $\omega = e^{2\pi i/3}$ and $\omega^2 = e^{4\pi i/3}$, which divide the circle into three equal arcs. We call this partition $\Pi_0$. For $k = 1, 2, \dots$, the partition $\Pi_k = f^{-k}(\Pi_0)$ divides the circle in $3 \cdot 2^k$ equal arcs.

We now define an analogous sequence of partitions of $\partial \Omega$. We define the {\em order} of an $\Omega$-horoball $\Omega_{p}$
as 
$$
\ord \Omega_p = \min_{v \in \partial \Omega_p} d_{\mathcal T_\infty}(v_{\Root}, v).$$
There are three $\Omega$-horoballs of order 0, which contain $v_{\Root}$. Inspection shows that for $k \ge 1$, there are $3 \cdot 2^{k-1}$ $\Omega$-horoballs of order $k$ and thus
$$3 + 3 + 6 + \dots + 3 \cdot 2^{k-1} = 3 \cdot 2^k$$
 $\Omega$-horoballs of order at most $k$. For $k =1, 2, \dots$, we define $\Lambda_k$ as the partition of $\partial \Omega$ into $3 \cdot 2^k$ arcs by the cusps of order $\le k$, i.e.~the points where $\Omega$-horoballs of order $\le k$ meet $\partial \Omega$. Since each arc in $\Lambda_k$ subtends the same number of edges of $\mathcal T_n$ up to an additive error of  $O(k)$, the harmonic measures of each arc in $\Lambda_k$ are equal and we have:

\begin{lemma}[Exterior structure lemma]
  \label{outer-structure}
  There is a  conformal mapping
  $$\psi: (\hat{\mathbb{C}} \setminus \overline{\mathbb{D}}, \infty) \to (\hat{\mathbb{C}} \setminus \overline{\Omega}, \infty)$$ which takes $\Pi_k$ to $\Lambda_k$ for any $k \ge 0$.
\end{lemma}

\subsection{Uniqueness of the Hausdorff limit}
\label{sec:uniqueness}

The interior and exterior structure lemmas (Lemmas \ref{inner-structure} and \ref{outer-structure}) show that any subsequential limit $\partial \Omega$ realizes the mating of $z \to \overline{z}^2$ and the Markov map $\rho(z)$ of the reflection group of an ideal triangle (see Section \ref{subsec:deltoid}).

The structure lemmas also show that Hausdorff limit of the $\mathcal T_n$ is unique.
Indeed, if $\mathcal T'_\infty \cup \partial \Omega'$ was another subsequential limit of $\mathcal T_n$, in addition to $\mathcal T_\infty \cup \partial \Omega$,
we could conformally map each complementary region in $\hat{\mathbb{C}} \setminus  (\mathcal T_\infty \cup \partial \Omega)$ to the corresponding complementary region in $\hat{\mathbb{C}} \setminus \mathcal (\mathcal T'_\infty \cup \partial \Omega')$. Lemmas \ref{inner-structure} and \ref{outer-structure} guarantee that these conformal mappings patch together to form a continuous self-map of the sphere $h: \hat{\mathbb{C}} \to  \hat{\mathbb{C}}$ which is conformal on $\hat{\mathbb{C}} \setminus (\mathcal T_\infty \cup \partial \Omega)$. The tree $\mathcal T_\infty$ is conformally removable as it is a union of real analytic arcs. Thus, $h$ is conformal on $\hat{\mathbb{C}} \setminus \partial \Omega$. By Lemmas \ref{weak-conformal-removability} and \ref{tree-square-sum}, $h$ must be a M\"obius transformation.

\subsection{Convergence of the Shabat polynomials}
\label{sec:convergence-shabat}

We  subdivide each $\Omega$-horoball $\Omega_{p_j}$ into triangles $\triangle(\vec{e}_i, p_j)$ by connecting the vertices of $\mathcal T_\infty$ on $\partial \Omega_{p_j}$ to $p_j$ by hyperbolic geodesics of $\Omega_{p_j}$. We colour the triangles $\triangle(\vec{e}_i, p_j) \subset \Omega$ black and white, so that adjacent triangles have opposite colours.

We conformally map
each black triangle $\triangle(\vec{e}_i, p)$ onto the upper half-plane $\mathbb{H}$ so that $\vec{e}_i \to [-1,1]$, $p_j \to \infty$ and each white triangle $\triangle(\vec{e}_i, p)$ onto the lower half-plane $\mathbb{L}$ so that $\vec{e}_i \to [-1,1]$, $p_j \to \infty$. Properties (F1) and (F2) from Section \ref{sec:relative}  guarantee that these
conformal maps glue together to form a holomorphic function $h$ defined on $\Omega$.
By choosing the colouring scheme appropriately, we can ensure that $h(v_{\Root}) = 1$ rather than $-1$.

From the description of the Shabat polynomials $p_n$ for the true trees $\mathcal T_n$ given in Section \ref{sec:shabat}, it is not difficult to see that the $p_n \to h$, uniformly on compact subsets of $\Omega$. Indeed, the Hausdorff convergence of $\mathcal T_n$ to $\mathcal T_{\infty}$ implies that the triangles $\triangle(\vec{e}_i, \infty)\subset \hat{\mathbb{C}}$ defined in Section \ref{sec:shabat} converge to the corresponding triangle $\triangle(\vec{e}_i, p)\subset \Omega_{p}$ in the Carath\'eodory topology. As $p_n$ and $h$ are conformal maps from these triangles to the upper or lower half-planes, this tells us that $p_n \to h$ uniformly on compact subsets of any triangle $\triangle(\vec{e}_i, p) \subset \Omega$. By considering a pair of triangles that have a common edge, one obtains that $p_n \to h$ uniformly on compact subsets of the union of these two triangles, which shows that $p_n \to h$ uniformly on compact subsets of $\Omega$ away from the vertices of the trees. Finally, one may use a similar argument to obtain the uniform convergence in a neighbourhood of each vertex $v \in \mathcal T_\infty$ by examining the behaviour of the maps $p_n$ and $h$ on the stars $\star_v$, which were defined in Section \ref{sec:trees-bounded-valence}.

\begin{figure}[h]\label{fig:h}
\centering
   \includegraphics[scale=0.18]{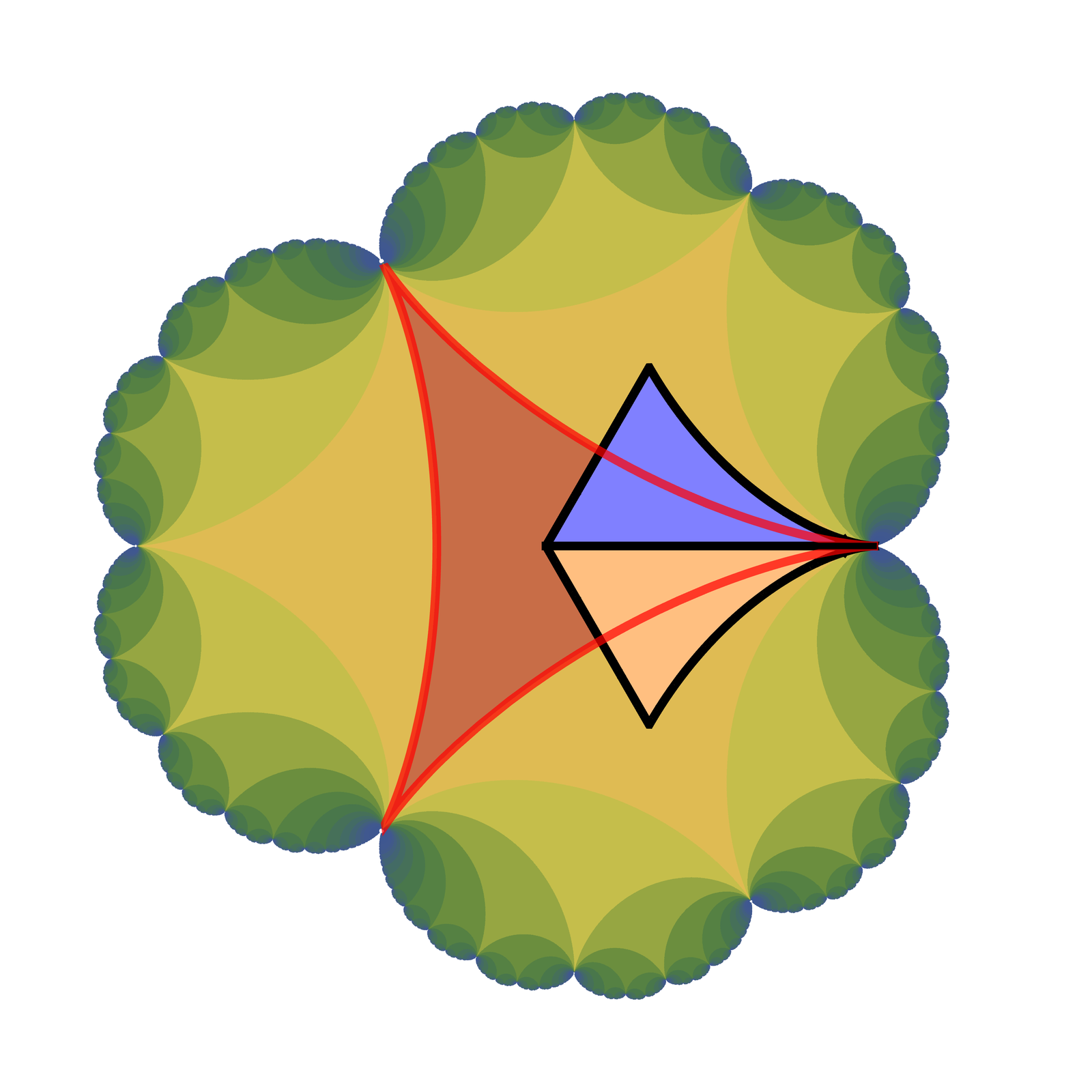} \qquad \includegraphics[scale=0.18]{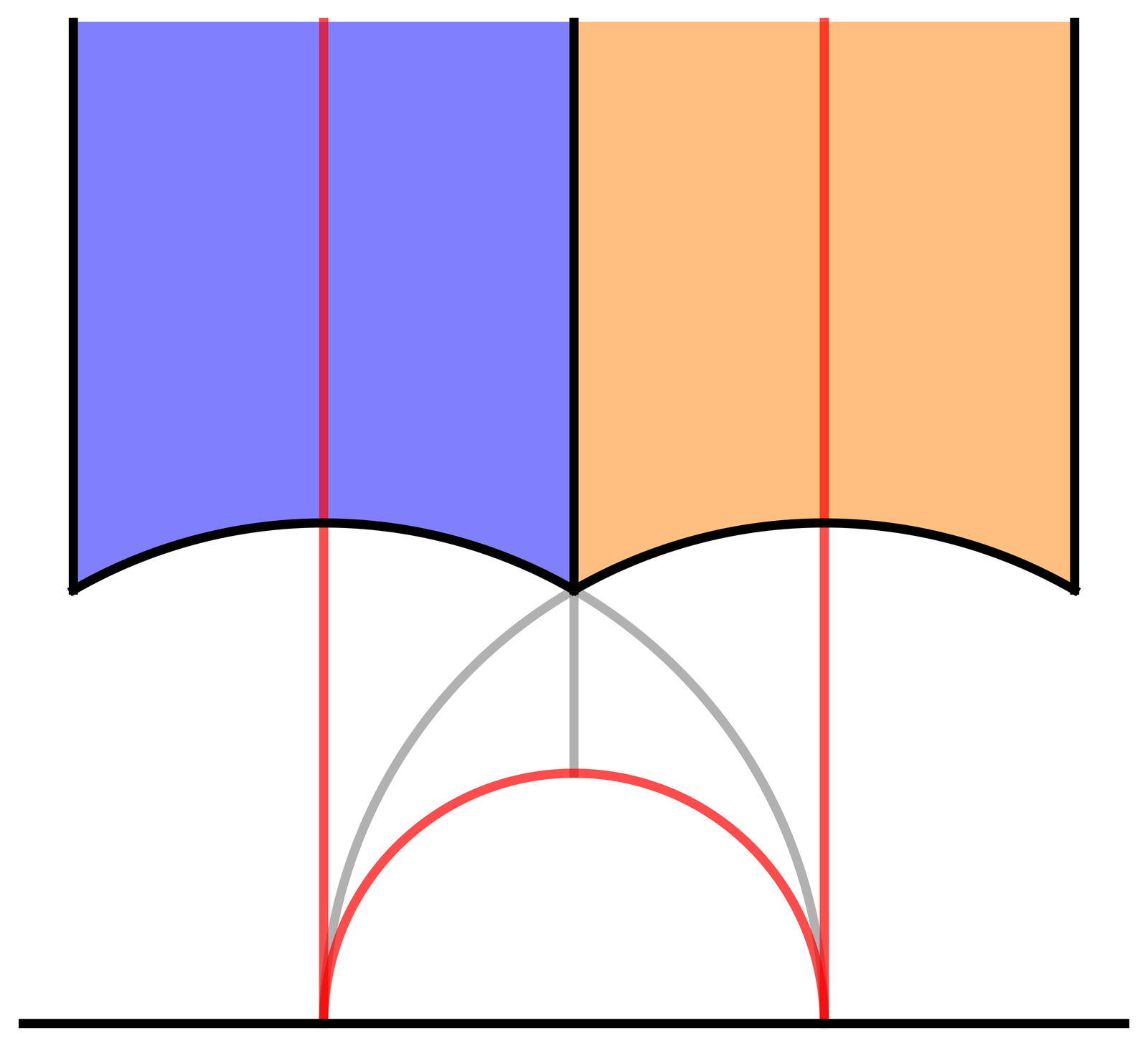}
    \caption{The fundamental domain for the function $h$, drawn in the upper half-plane instead of the disk, is twice as big as the fundamental domain for $\mathbb{H}/\PSL(2, \mathbb{Z})$. The blue part of the fundamental domain is sent to the upper half-plane $\mathbb{H}$, while the orange part is send to the lower half-plane $\mathbb{L}$.}
\end{figure}

If $R$ is the Riemann map from $\mathbb{D} \to \Omega$, then $h \circ R$ is a function on the unit disk whose fundamental domain consists of two copies of the fundamental domain for
$\PSL(2, \mathbb{Z})$, see Figure \ref{fig:h}.

\appendix

\section{Trivalent true trees are dense}
\label{sec:true-trees-are-dense}

Let $K$ be a connected compact set in the plane. In this appendix, we show that one can approximate $K$ in the Hausdorff topology by finite trivalent true trees, thereby giving another proof of Bishop's theorem.

\begin{figure}[h]
\centering
   \includegraphics[scale=0.24]{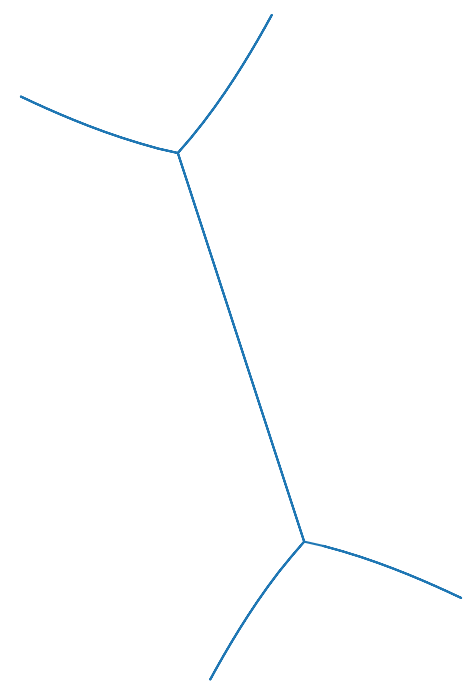} \quad \includegraphics[scale=0.18]{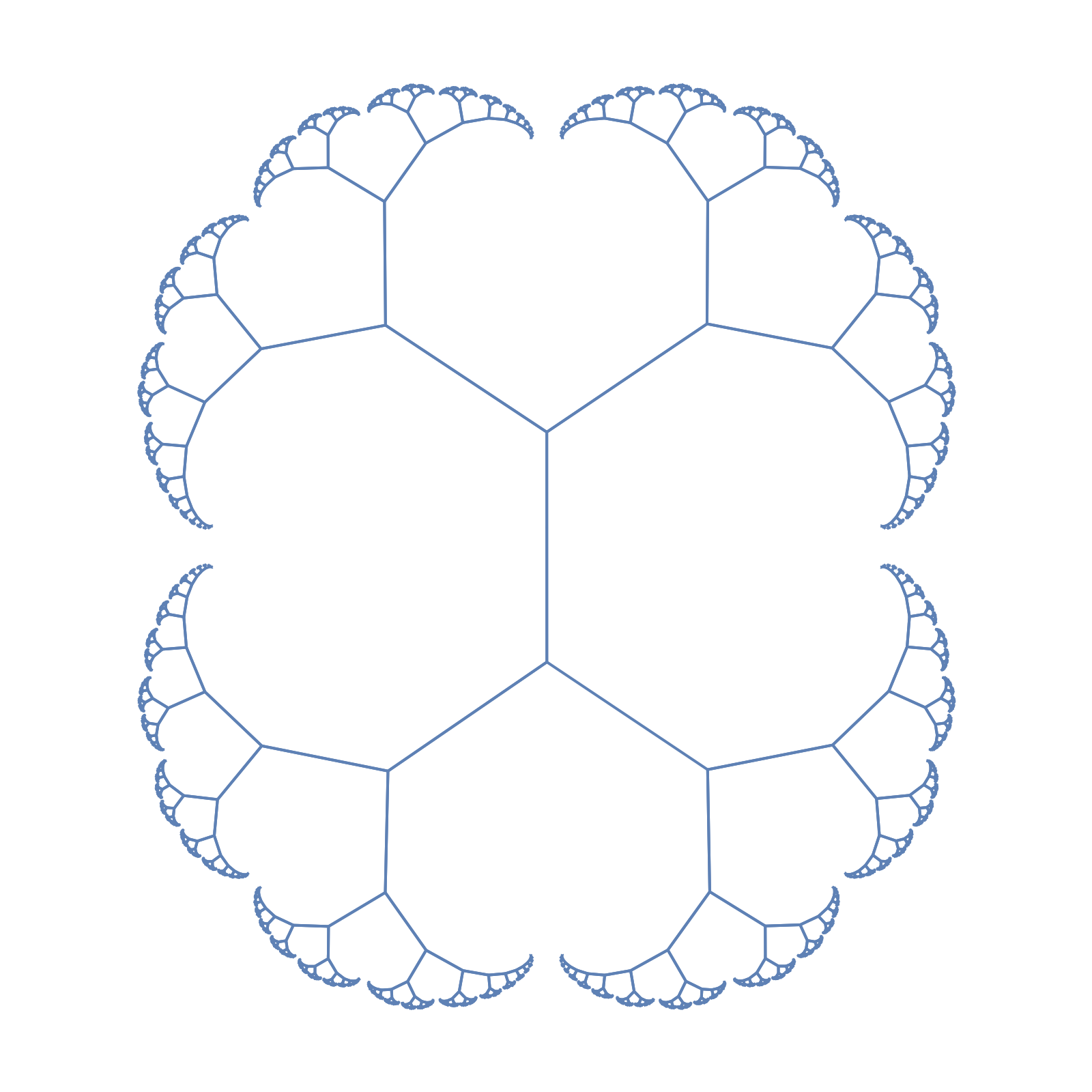}
 \caption{Unbalanced truncations of the infinite trivalent tree.}
 \label{fig:fake-deltoid}
\end{figure}

 Start with a finite trivalent tree $\mathcal T'_1$, for instance, with the tree on the left side of Fig.~\ref{fig:fake-deltoid} which consists of five edges.
At each step, add two edges to each boundary vertex. This gives us a sequence of true trees $\{\mathcal T'_n\}_{n=1}^\infty.$
The arguments presented in this paper show that the finite trees $\mathcal T'_n$ converge in the Hausdorff topology to an infinite trivalent tree union a Jordan curve: $\mathcal T'_\infty \cup \partial \Omega'$.

\begin{figure}[h]
\centering
   \includegraphics[scale=0.6]{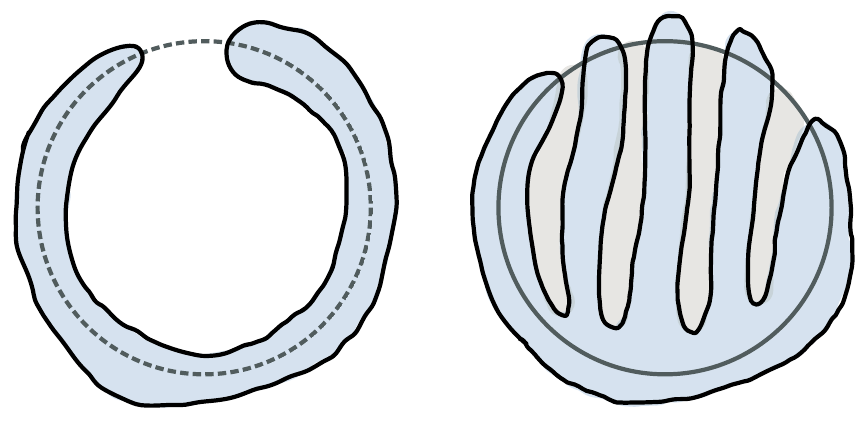}
 \caption{Approximating the unit circle and the unit disk in the Hausdorff topology by thin Jordan domains.}
 \label{fig:hausdorff-approximation}
\end{figure}

Our objective is to show that for any compact connected set $K \subset \mathbb{C}$ and $\varepsilon > 0$, one can choose the starting tree $\mathcal T'_1$ appropriately so that the Hausdorff distance $$d_{\Haus}(L \circ \mathcal T_n', K) < \varepsilon,$$ for some linear mapping $L(z) = az+b$ in $\aut \mathbb{C}$. (The linear mapping compensates for the fact that the conformal map to $\mathbb{C} \setminus \mathcal T'_n$ is hydrodynamically normalized.) We now make two reductions.

\medskip

{\em Reduction 1.} It is enough to show that for any Jordan curve $\gamma$, one can choose the starting tree $\mathcal T'_1$ so that $d_{\Haus}(L \circ \partial \Omega', \gamma) < \varepsilon$.

\medskip

Indeed, one can approximate any compact connected set $K \subset \mathbb{C}$ in the Hausdorff topology by Jordan curves $\gamma_k$ that are $(1/k)$-thin, i.e.~any point in the domain $\Gamma_k$ enclosed by $\gamma_k$ lies within $1/k$ of $\gamma_k$. This ensures that the domains $\Gamma_k$ converge in the Hausdorff topology to $K$. See Fig.~\ref{fig:hausdorff-approximation} above for examples.

Therefore, if $\mathcal T'_{k,l}$ is a sequence of infinite trees such that $d_H( L_{k,l} \circ \partial \Omega'_{k,l}, \gamma_k) \to 0$, then
 $d_H(L_{k,l} \circ \mathcal T'_{k,l}, \gamma_k) \le 2/k$ for all sufficiently large $l \ge l_0(k)$. A diagonal argument produces a sequence
 of infinite trees which converges to $K$ after linear rescaling.

\medskip

{\em Reduction 2.} We may assume that $K = \partial \tilde{\Omega}$ is the image of the boundary of the developed deltoid $\partial \Omega$ under a quasiconformal map $f: \mathbb{C} \to \mathbb{C}$, which is conformal on $\Omega$.

\medskip

It is well known that quasiconformal images of the unit circle (quasicircles) are dense in the collection of Jordan curves. It turns out that quasiconformal images of any fixed Jordan curve are also dense. An anti-symmetrization argument from \cite{kuhnau, smirnov} allows one to choose $f$ to be conformal on $\Omega$. See Lemma \ref{quasicircles-are-dense1} below.

\medskip

\subsection{Trivalent tree weldings}

Recall that any non-root vertex in the trees $\mathcal T_n$ and $\mathcal T_\infty$ can be labeled by a digit $1, 2, 3$ followed by a sequence of left and right turns. In order to label the vertices of $\mathcal T'_n$ and $\mathcal T'_\infty$ in a similar fashion, we designate a vertex in $\mathcal T'_1$ as the root vertex and select one of the adjacent vertices as the vertex labeled 1.

Let $\varphi: (\mathbb{D}, 0, 1) \to (\Omega, v_{\Root}, p_1)$ and $\psi: (\mathbb{D}_e, \infty, 1) \to (\Omega_e, \infty, p_1)$ be conformal mappings to the interior and exterior of the developed deltoid respectively, where $$p_1 \, = \, \lim_{k \to \infty} v_{1R^k} \, = \, \lim_{k \to \infty} v_{3L^k}$$ is one of the three cusps of the developed deltoid of order 0. The composition $h = \psi^{-1} \circ \varphi: \partial \mathbb{D} \to \partial \mathbb{D}$ defines a homeomorphism of the unit circle, which is called the {\em welding homeomorphism} of $(\partial \Omega, v_{\Root}, p_1, \infty)$. Form the analogous mappings $\varphi', \psi'$ and $h'$ for $(\partial \Omega', v'_{\Root}, p'_1, \infty)$. Inspection shows that the weldings $h$ and $h'$ are related by a piecewise linear homeomorphism $F$ of the unit circle: $h' = F \circ h$.

For instance, in the example depicted in Fig.~\ref{fig:fake-deltoid}, to describe $F$, we divide the unit circle into three equal arcs and map these onto arcs of lengths $\pi, \pi/2, \pi/2$ respectively, which indicates the fact that one third of the tree has the same number of edges as the other two thirds. Let $\TPL$ denote the collection of piecewise linear homeomorphisms of the unit circle that arise in this way and $\TW = \{ F \circ h : F \in \TPL \}$ be the collection of all trivalent tree weldings.

It is not difficult to see that for any quasisymmetric homeomorphism of the unit circle
$F \in \QS_1$ which fixes $1 \in \partial \mathbb{D}$, there is a sequence of homeomorphisms $F_k \in \TPL$ whose quasisymmetry constants are uniformly bounded such that $$F_k^{-1} \circ F: \partial \mathbb{D} \to \partial \mathbb{D}$$ tend uniformly to the identity. (One may choose the  homeomorphisms $F_k$ so that their quasisymmetry constants are comparable to the quasisymmetry constant of $F$.)

\subsection{An overview of the proof}

The curve $K = \partial \tilde{\Omega}$ divides the Riemann sphere $\sph$ into an interior domain $\tilde{\Omega}$ and an exterior domain
$\tilde{\Omega}_e$.
Composing $f$ with conformal maps $\psi: (\mathbb{D}_e, \infty, 1) \to (\Omega_e, \infty, p_1)$ and $\tilde{\psi}^{-1}: (\tilde{\Omega}_e, \infty, f(p_1)) \to (\mathbb{D}_e, \infty, 1)$, we get a quasiconformal self-mapping of the exterior unit disk
$$
F = \tilde{\psi}^{-1} \circ f \circ \psi.
$$
The following lemma (whose proof will be presented in Section \ref{sec:PL-strip}) allows us to approximate $F$ by quasiconformal self-maps $F_k$ of $\mathbb{D}_e$ with
$F_k|_{\partial \mathbb{D}} \in \TPL$:

\begin{lemma}
\label{PL-approximation}
Let $F$ be a quasiconformal self-map of the exterior of the unit disk which fixes $1$ and $\infty$.
We can approximate $F$ by quasiconformal maps $F_k: (\mathbb{D}_e, 1, \infty) \to (\mathbb{D}_e, 1, \infty)$ so that:
\begin{enumerate}
\item $F_k = F$ on $\{z > 1 + 1/k \}$.
\item The dilatations $\| \mu_{F_k} \|_\infty < c < 1$ are uniformly bounded.
\item Restricted to the unit circle, $F_k$ is one of the piecewise linear homeomorphisms described above that relates the welding of the ``genuine'' developed deltoid $\Omega$ and a ``generalized'' developed deltoid $\Omega_k'$.
\end{enumerate}
\end{lemma}

By Properties 1 and 2, we can select quasiconformal maps $f_k: \mathbb{C} \to \mathbb{C}$ which tend to $f$, are conformal on $\Omega$ and
have dilatations  $\psi_* \mu_{F_k}$ on $\Omega_e$. (Since $\partial \Omega$ has area zero by Lemma \ref{tree-square-sum}, the quasiconformal map $f_k$ is uniquely specified by its dilation off $\partial \Omega$ up to post-composition with a M\"obius transformation.)
We set $\partial \tilde{\Omega}_k = f_k(\partial \Omega)$. Property 3 tells us that  $\bigl (\partial {\tilde{\Omega}}_k, f_k(v_{\Root}), f_k(p_1), \infty \bigr )$ has welding homeomorphism $h_k = F_k \circ h \in \TW$. 

Let $\Omega'_k$ be the generalized developed deltoid with welding $h_k \in \TW$.
We now use partial conformal removability techniques to show that $\tilde{\Omega}_k = L(\Omega'_k)$ for some linear map $L \in \aut \mathbb{C}$.
 Since  $\partial \tilde{\Omega}_k$ and $\partial \Omega'_k$ realize the same welding homeomorphism, the conformal map 
 $$
 \phi: \bigl (\sph \setminus \partial \tilde{\Omega}_k, f_k(v_{\Root}), f_k(p_1), \infty \bigr ) \to  \bigl ( \sph \setminus \partial \Omega'_k, v'_{\Root;\;\! k}, p'_{1;\;\!k}, \infty \bigr )
 $$
  extends continuously to a homeomorphism on the Riemann sphere. We claim that $\phi$ is conformal on the whole Riemann sphere (and therefore, a linear mapping since it fixes infinity).

As in Section \ref{sec:uniqueness-estimate}, for each non-root vertex $v$ of $\mathcal T'_\infty \cong \mathcal T_\infty$, one can define the shadow
$s'_v \subset \partial \Omega'$ as the shorter arc of $\partial \Omega'$ between $vLRL^\infty$ and $vRLR^\infty$. Then,
each point of $\partial \Omega'$ which is not a cusp lies in infinitely many shadows and
\begin{equation}
\label{eq:square-sum1}
\sum_{v \ne v_{\Root}} \diam^2 s'_v < \infty.
\end{equation}

We can also define a collection of shadows $\{\tilde{s}_v\}$ for $\tilde{\Omega}_k$ by taking the images of the shadows for the developed deltoid $\Omega$ under the quasiconformal map $f_k$. Since in each case, the shadows are defined in terms of the combinatorics of the infinite trivalent tree,
 $\phi$ takes $\tilde{s}_v$ onto $s'_v$ for each vertex $v \ne v_{\Root}$.
Recall that in Section \ref{sec:uniqueness-estimate}, to each shadow $s_v \in \partial \Omega$, $v \ne v_{\Root}$, we associated a round set $B_v$
with
$$
\area B_v \asymp \diam^2 B_v, \qquad \diam B_v \asymp \diam s_v, \qquad \dist(B_v, s_v) \lesssim \diam s_v,
$$
so that $\{B_v\}$ are disjoint and contained in a bounded set. From the quasisymmetry of $f_k$, we deduce that
\begin{equation}
\label{eq:square-sum2}
\sum_{v \ne v_{\Root}} \diam^2 \tilde{s}_v < \infty.
\end{equation}

In view of Lemma \ref{weak-conformal-removability}, the inequalities (\ref{eq:square-sum1}) and (\ref{eq:square-sum2}) show that $\phi$ is a M\"obius transformation. This completes the proof of Bishop's theorem, modulo the technical Lemmas  \ref{PL-approximation} and \ref{quasicircles-are-dense1} which will be proved below.

\subsection{Quasiconformal images of Jordan curves}
\label{sec:qc-jordan}

In the following lemma, we explain how to approximate Jordan curves by quasiconformal images of a given Jordan curve:

\begin{lemma}
\label{quasicircles-are-dense1}
Let $\Omega \subset \mathbb{C}$ be a bounded Jordan domain.
For any Jordan curve $\gamma$ and $\varepsilon > 0$, one can find a quasiconformal map $f: \mathbb{C} \to \mathbb{C}$, which is conformal on $\Omega$ and takes
 $\partial \Omega$ onto a Jordan curve $\partial \tilde{\Omega}$ for which the Hausdorff distance $d_H(\partial \tilde{\Omega}, \gamma) < \varepsilon$.
 \end{lemma}

Before proving the above lemma, we first make a preliminary observation:

\begin{lemma}
\label{quasicircles-are-dense2}
Let $\gamma \subset \mathbb{C}$ be a smooth Jordan curve. For any $\delta, \varepsilon > 0$, one can express $\gamma$ as the image of the unit circle under a quasiconformal map $f: \mathbb{C} \to \mathbb{C}$ such that $f(A(0; 1-\delta, 1+\delta))$ contains an $\varepsilon$-neighbourhood of $\gamma$.
\end{lemma}

 \begin{proof}
It is not difficult to express $\gamma$ as the image of the unit circle under a smooth quasiconformal map $f_1: \mathbb{C} \to \mathbb{C}$. Pick $\rho > 0$
so that $f_1(A(0; \rho, 1/\rho))$ contains an $\varepsilon$-neighbourhood of $\gamma$.
The desired quasiconformal map can be obtained by a radial reparametrization of $f_1$:
$$
f(re^{i\theta}) := f_1(\phi(r) e^{i\theta}),
$$
where
$\phi: [0, \infty) \to [0, \infty)$ is a homeomorphism which is identity outside $[\rho, 1/\rho]$ and takes $[\rho, 1/\rho]$ to $[1-\delta, 1+\delta]$.
 \end{proof}

 \begin{proof}[Proof of Lemma \ref{quasicircles-are-dense1}]
Since smooth Jordan curves are dense in the set of all Jordan curves, one can find two smooth quasiconformal maps
 $f_1, f_2: \mathbb{C} \to \mathbb{C}$ such that $d_H(f_1(\partial \mathbb{D}), \partial \Omega) < \varepsilon/2$ and $d_H(f_2(\partial \mathbb{D}), \gamma) < \varepsilon/2$. Anti-symmetrizing as in \cite{kuhnau, smirnov}, one may assume that $f_1, f_2$ are conformal on the unit disk. The idea is to take $f = f_2 \circ f_1^{-1}$ and $\partial \tilde{\Omega} =  f_2 \circ f_1^{-1}(\partial \mathbb{D})$.
 
 To make this argument work, one has to be slightly careful when choosing the maps $f_1$ and $f_2$. 
 Here is the precise construction:
 \begin{enumerate}
\item We first choose $f_2$ so that $d_H(f_2(\partial \mathbb{D}), \gamma) < \varepsilon/2$.

\item We then choose $\delta > 0$ sufficiently small so that $f_2(A(0; 1-\delta, 1+\delta))$ is contained in an $\varepsilon/2$-neighbourhood of $f_2(\partial \mathbb{D})$, and thus in an $\varepsilon$-neighbourhood of $\gamma$.

\item Finally, we use Lemma \ref{quasicircles-are-dense2} to select a quasiconformal map $f_1$ so that 
$$f_1^{-1}(\partial \Omega) \subset A(0; 1-\delta, 1+\delta).$$
 \end{enumerate}
It is clear from the construction that the composition $f = f_2 \circ f_1^{-1}$ maps $\partial \Omega$ into the $\varepsilon$-neighbourhood of $\gamma$.
 \end{proof}
 
\subsection{Piecewise-linear approximations} 
  \label{sec:PL-strip}
  
  In the following lemma, we explain how to extend quasisymmetric maps from the boundary of a horizontal strip to the interior. This is a special case of a result of V\"{a}is\"{a}l\"{a}, see \cite{vaisala}.
  
 \begin{lemma}
 \label{PL-strip}
Let   $\mathcal S = \{ (x,y) \in \mathbb{R}^2 : 0 < y < 1 \}$ be a horizontal strip of width 1.
 Suppose $\phi_0, \phi_1: \mathbb{R} \to \mathbb{R}$ are  $k$-quasisymmetric maps which move points a bounded distance, i.e.~$|\phi_i(x) - x| < C$, $i=0,1$. There exists a $k_1$-quasiconformal map $\Phi: \mathcal S \to \mathcal S$ which takes $(x,0) \to (\phi_1(x),0)$ and $(x,1) \to (\phi_2(x), 1)$, with $k_1$ depending only on $k$ and $C$.
 \end{lemma}
 
 \begin{proof}
 We begin by defining $\Phi$ on $\partial \mathcal S$ by  $(x,0) \to (\phi_0(x),0)$ and $(x,1) \to (\phi_1(x), 1)$.
We can partition  $\mathcal S$ into a union of squares $\{S_n\}$ using the vertical segments $\{\ell_n\}_{n \in \mathbb{Z}}$ which connect $(n, 0)$ and $(n,1)$. Similarly, we can partition $\mathcal S$ into a union of conformal rectangles $\{\tilde{S}_n\}$ using the line segments $\{\tilde{\ell}_n\}_{n \in \mathbb{Z}}$  which connect $(\phi_0(n), 0)$ to $(\phi_1(n), 1)$. We first extend $\Phi$ to the vertical segments  $\{\ell_n\}_{n \in \mathbb{Z}}$, so that it is linear on each segment $\ell_n$ and takes $\ell_n$ to $\tilde{\ell}_n$. For each $n \in \mathbb{Z}$, we extend $\Phi$ from $\partial S_n \to \partial \tilde{S}_{n}$ to $\Phi: S_n \to \tilde{S}_{n}$ using the Beurling-Ahlfors extension.
(The assumption on the maps $\phi_0$ and $\phi_1$ guarantees that $\partial \tilde{S}_n$ are uniform quasicircles.)
 \end{proof}
 
 We now show how to approximate quasiconformal self-maps of $\mathbb{D}_e$ by ones that are piecewise-linear on the unit circle. In the proof below, we will use a variant of the above lemma for the annulus $A(0; 1, 1+1/k)$\,:
 
 \begin{proof}[Proof of Lemma \ref{PL-approximation}]
The idea is to define $F_k = F \circ \Phi_k$ by composing $F$ with a quasiconformal homeomorphism $\Phi_k: \mathbb{D}_e \to \mathbb{D}_e$ which is identity on $|z| > 1+1/k$.

Let $\Lambda_k \in \TPL$ be a piecewise linear map whose quasisymmetric constant is comparable to that of $F$ such that
 $\phi_k = F^{-1} \circ \Lambda_k$ moves points on the unit circle by $O(1/k)$. By Lemma \ref{PL-strip}, $\phi_k$ admits a quasiconformal extension $\Phi_k$ to the exterior unit disk which is identity on $\{z > 1 + 1/k\}$.
From the construction, it is clear that $F_k = F \circ \Phi_k$ satisfies Properties 1--3 as desired. 
\end{proof}
 
\section{True tree approximation of cauliflower}

In this appendix, we describe a sequence of true trees, whose limit is the cauliflower, the Julia set of $f(z) = z^2+1/4$.
Since the arguments are similar to the ones for the finite truncations of the infinite trivalent tree, we only give a brief sketch of the proofs, with an emphasis on the differences.

 \begin{figure}[h!]
\centering
\includegraphics[scale=0.65]{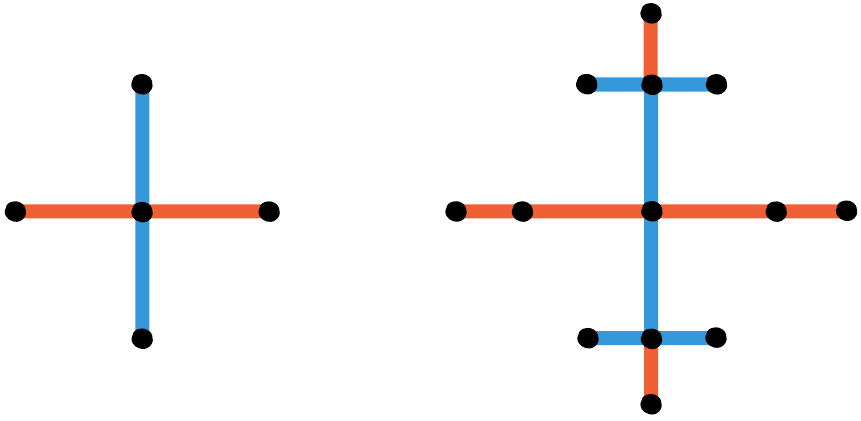}
\caption{A sequence of true trees given by an inductive construction.}
\label{fig:cauliflower-trees2}
 \end{figure}
 
Let $T_1$ be a planar tree  which consists of a root vertex $v_{\Root}$ and four edges
$$
\overline{v_{\Root}v_\uparrow}, \quad
 \overline{v_{\Root}v_{\rightarrow}}, \quad
 \overline{v_{\Root}v_{\downarrow}}, \quad
 \overline{v_{\Root}v_{\rightarrow}},
 $$ labeled counter-clockwise. We colour the edges  $\overline{v_{\Root}v_{\uparrow}}, \  \overline{v_{\Root}v_{\downarrow}}$ blue and $\overline{v_{\Root}v_{\leftarrow}}, \ \overline{v_{\Root}v_{\rightarrow}}$ red.
 To form $T_{n+1}$ from $T_n$, we attach additional edges at each leaf vertex:
\begin{itemize}
\item If a leaf edge is red, we attach another red edge at the leaf vertex.
\item  If a leaf edge is blue, we attach three edges, coloured blue-red-blue in counter-clockwise order.
 \end{itemize}

The trees $T_1$ and $T_2$ are depicted on Fig.~\ref{fig:cauliflower-trees2}. 
 From this description, it is easy to see that $T_n$ is made out of
 $$
 4 + 8 + 16 + \dots + 2^{n+1} = 2^{n+2} - 4
 $$
  edges, with the same number of red and blue edges. 
 
 Let $\mathcal T_n$ be the true tree representative of $T_n$. Note that the colouring is only used to describe the combinatorics of $T_n$, it plays no role in how the true tree $\mathcal T_n$ is constructed from $T_n$.

 \begin{figure}[h!]
\centering
\includegraphics[scale=0.25]{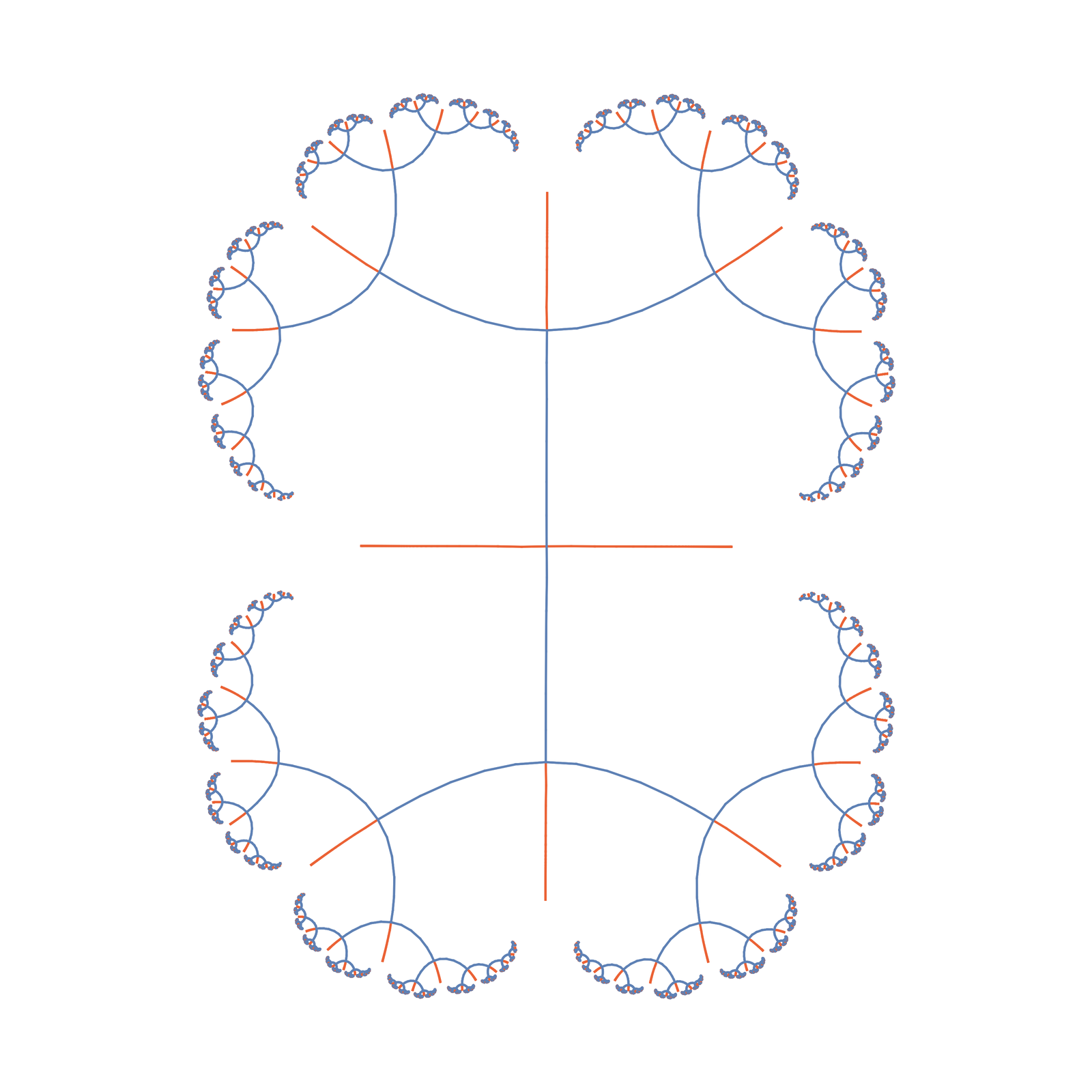}
\includegraphics[scale=0.25]{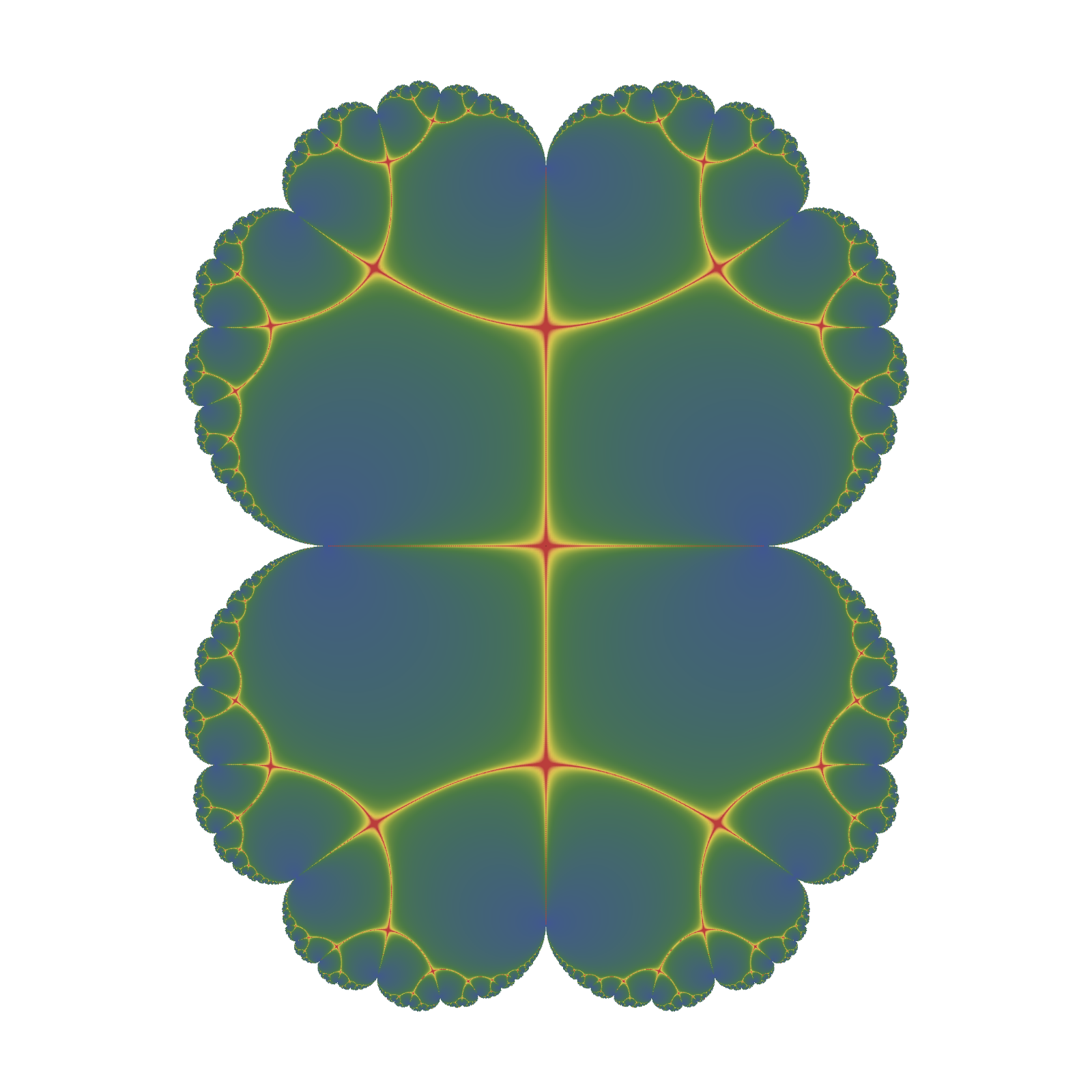}
\caption{A sequence of true trees which approximates $\mathcal J(z^2+1/4)$.}
\label{fig:cauliflower-approximation}
 \end{figure}

\begin{theorem}
\label{main-thm2}
The trees $\mathcal T_n$ converge in the Hausdorff topoology to an infinite tree union a Jordan curve
$\mathcal T \cup \partial \Omega$. The Jordan curve $\partial \Omega$ is the Julia set of $z^2+1/4$, while the set of vertices of $\mathcal T$ is the grand orbit of the critical point 0 of $f(z) = z^2+1/4$. Let $\psi: \Omega \to \mathbb{C}$ be the Fatou coordinate at the parabolic fixed point $1/2 \in \mathcal J(f)$, with 
$$
\psi(f(z)) = \psi(z) + 1, \qquad \psi(0) = 0.
$$
The Shabat polynomials $p_n(z)$ of $\mathcal T_n$, with $p_n(0) = 1$, converge uniformly on compact subsets of $\Omega$ to $\cos(\pi \cdot \psi(z))$.
\end{theorem}

\subsection{Topology of a subsequential limit}
\label{sec:july}

We first show that any Hausdorff limit of the trees $\mathcal T_n$ is ambiently homeomorphic to the set depicted on the right side of Fig.~\ref{fig:cauliflower-approximation}. Since each vertex of $\mathcal T_n$ has at most 4 neighbours and the sums
\begin{equation}
\label{eq:square-sum-cauliflower}
S_n = \sum_{e \in \mathcal T_n} s(e)^2
\end{equation}
are uniformly bounded above, we are in the setting of Theorem \ref{edges-do-not-shrink-thm}. For any blue edge $e_0$, the numbers $s(e_0)$ are uniformly bounded below, so the blue edges do not shrink.
By Lemma \ref{basic-trees2}, the red edges also do not shrink. 
Therefore, any subsequential Hausdorff limit of $\mathcal T_n$ contains an infinite tree $\mathcal T_\infty$ whose edges are real-analytic arcs.

Let $\mathcal B_n \subset \mathcal T_n$ be the subtree consisting of blue edges. Arguing as in Section \ref{sec:shrinking-diameters}, one can show that the subsequential limit of the subtrees $\mathcal B_n$  is an infinite tree union a Jordan curve $\partial \Omega$. 

Perhaps the new feature of the trees $\mathcal T_n$ are the red edges, so we discuss their behaviour in more detail.
The red edges are naturally grouped into {\em twigs}. There are two twigs emanating from the root vertex, which we denote $\twig_{v_{\Root}, \leftarrow}$ and $\twig_{v_{\Root}, \rightarrow}$, while 
a single twig $\twig_v$ emanates from each degree 4 vertex $v$, other than the root vertex, which we denote by $\twig_v$.

The following lemma says that each twig connects a degree 4 vertex in $\mathcal T_n$ to a cusp in $\partial \Omega$, where it meets the two enclosing blue branches, without protruding outside of $\Omega$\,:

\begin{lemma}
\label{twig-lemma}
Let $\twig_v^{(n)}$ be a twig in $\mathcal T_n$.

{\em (i)} Any Hausdorff limit of $\twig_v^{(n)}$ is contained in $\overline{\Omega}$. 

{\em (ii)} Any Hausdorff limit of $\twig_v^{(n)}$ connects $v \in \mathcal T_\infty$ to the cusp $p_v \in \partial \Omega$. 
\end{lemma}
 
 \begin{figure}[h!]
\centering
\includegraphics[scale=0.33]{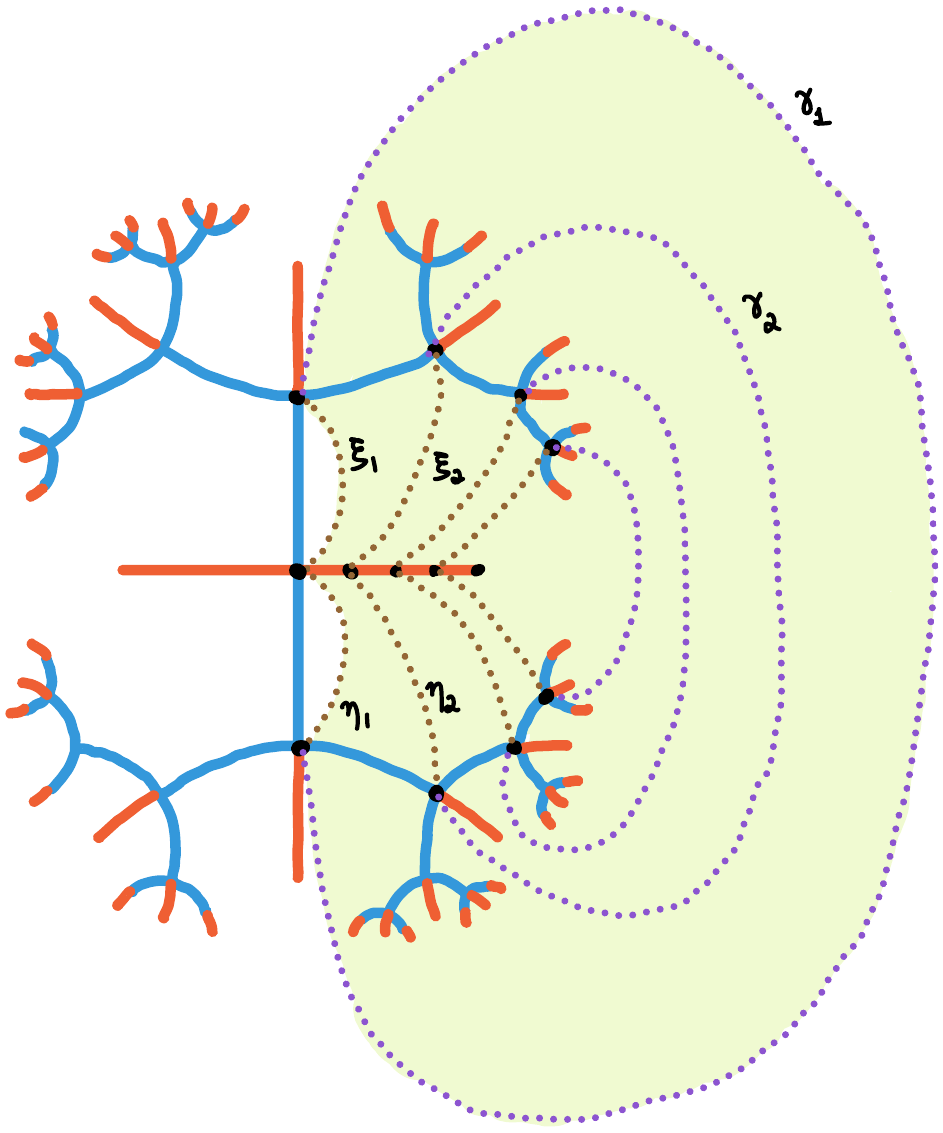}
\caption{The twigs are enclosed by the hyperbolic geodesics $\gamma_m^{(n)}$.}
\label{fig:enclosing-geodesics}
 \end{figure}

\begin{proof}[Sketch of proof]
We explain the argument for the twig $\twig_{v_{\Root}, \rightarrow}$ as the general case is similar. We pass to a subsequence so that $\mathcal B_n$ and
 $\twig_{v_{\Root}, \rightarrow}$ converge in the Hausdorff topology as $n \to \infty$.  We denote the associated cusp by $$p_{v_{\Root}, \rightarrow} = \lim_{m \to \infty} (\uparrow\! R^{m-1}) 
 = \lim_{m \to \infty} (\downarrow\! L^{m-1}) \in \partial \Omega.$$
 
(i) For $1 \le m \le  {n-1}$, we construct hyperbolic geodesics $\gamma_m^{(n)} \subset \hat{\mathbb{C}} \setminus \mathcal T_n$ connecting
$\uparrow\! R^{m-1}$ and $\downarrow\! L^{m-1}$ as in Figure~\ref{fig:enclosing-geodesics}.
By construction, $\twig_{v_{\Root}, \rightarrow}^{(n)} \subset \mathcal T_n$ is contained in the subdomain of $\hat{\mathbb{C}} \setminus \mathcal B_n$ enclosed by $\gamma_m^{(n)}$.

For any $m \ge 1$, the Hausdorff limit of the geodesics $\gamma_m^{(n)}$ as $n \to \infty$ is composed of three pieces: two pieces $\gamma_{m,1}, \gamma_{m, 3}$ are hyperbolic geodesics in the tiles that make up $\Omega$, while the middle piece is
a  hyperbolic geodesic $\gamma_{m,2}$ in $\hat{\mathbb{C}} \setminus \overline{\Omega}$ which connects two cusps $p_{-m}, p_m \in \partial \Omega$.
Therefore, the Hausdorff limit of the twigs $\twig_{v_{\Root}, \rightarrow}^{(n)}$ is contained in $\overline{\Omega}$ union the subdomain of $\hat{\mathbb{C}} \setminus \overline{\Omega}$ enclosed by $\gamma_{m,2}$.

Since the points $p_{-m}, p_m$ tend to the cusp $p_{v_{\Root}, \rightarrow} \in \partial \Omega$ as $m \to \infty$, the subdomains  of $\hat{\mathbb{C}} \setminus \overline{\Omega}$ enclosed by $\gamma_{m,2}$ shrink down to $p_{v_{\Root}, \rightarrow}$. It follows that the Hausdorff limit of the twigs $\twig_{v_{\Root}, \rightarrow}^{(n)}$ is contained in $\overline{\Omega}$ as desired.

(ii) Fig.~\ref{fig:enclosing-geodesics} depicts a decreasing sequence of simply-connected domains
$$
W_1^{(n)} \,  \supset \, W_2^{(n)} \, \supset \, \dots \, \supset \, W_{m-1}^{(n)},
$$
which contain the set $\mathcal T_n(v_{\rightarrow^m}) \cup \bigl \{ \uparrow\! R^{m-1}, {\downarrow\! L^{m-1}} \bigr \}$.
A moduli estimate similar to the one in Section \ref{sec:shrinking-diameters} shows that 
$$
\Mod \bigl (W_j^{(n)} \setminus W_{j+1}^{(n)} \bigr ) \gtrsim 1/j, \qquad j = 1, 2, \dots, m-2.
$$
By the parallel rule,
$$
\Mod \bigl ( W_1^{(n)} \setminus W_{m-1}^{(n)} \bigr )
\, \gtrsim \,
1 + 1/2  + \dots + 1/(m-2) \, \asymp \, \log m.
$$
Since the initial domain  $W_1^{(n)}$ is contained in a ball $B(0,R_0)$ where $R_0 > 0$ is a universal constant, Lemma \ref{lem:egg_yolk} implies that the diameter of $\mathcal T_n(v_{\rightarrow^m}) \cup \bigl \{ \uparrow\! R^{m-1}, {\downarrow\! L^{m-1}} \bigr \}$ is small, which means that the red twig and the two blue branches come together at $p_{v_{\Root}, \rightarrow} \in \partial \Omega$.
\end{proof}

\subsection{Tile decomposition}

As shown on the right side of Fig.~\ref{fig:cauliflower-approximation}, the repeated pre-images of the line segment $[-1/2,1/2]$ separate $\Omega$, the interior of the filled Julia set of $f(z) = z^2+1/4$, into a countable collection of {\em tiles}.  The union of  these curves forms a tree whose vertices are points in the grand orbit of the critical point 0. We designate the critical point 0 as the root vertex.
Note that $[0,1/2]$ is not a single edge but the union of countably many edges:
$$
[0,1/2] \, = \, [0, f(0)] \, \cup \, [f(0), f^{\circ 2}(0)] \, \cup \, [f^{\circ 2}(0), f^{\circ 3}(0)] \, \cup \, \dots
$$
We label the tiles as $\Omega_{p, L}$ or $\Omega_{p, R}$, where $p$ ranges over the cusps in $\partial \Omega$.
A {\em bi-tile} $\Omega_p$ is a horoball-like region formed by taking the interior of the closure of $\Omega_{p,L} \cup \Omega_{p,R}$. Thus, $\Omega$ is organized into a union of bi-tiles, as well as a union of tiles.

 Under iteration, any tile is eventually mapped onto $\Omega_{1/2,L}$ or $\Omega_{1/2,R}$. The tiles $\Omega_{1/2,L}$ or $\Omega_{1/2,R}$ are invariant under $f$, and $f$ restricts as a conformal automorphism on $\Omega_{1/2,L}$ and $\Omega_{1/2,R}$. We record the following two properties of $\Omega$, which come from the dynamics of $f$ and the symmetry of $\Omega$ with respect to the real axis:
\begin{enumerate}[leftmargin=1.75cm, label={\rm (CT\arabic*)}]
\item \label{item:cauliflower-tree1} If $\Omega_{p,X}$ is a tile, then each edge in $\partial \Omega_{p,X}$ has the same relative harmonic measure as viewed from $p$, i.e.~if $e_1, e_2 \subset \Omega_{p,X}$, then
$$
\lim_{z \to p, \, z \in \Omega_{p,X}} \frac{\omega_z(e_1)}{\omega_z(e_2)} = 1.
$$
\item \label{item:cauliflower-tree2} If $e$ is an edge that belongs to two neighbouring tiles $\Omega_{p,X}$ and $\Omega_{q,Y}$, then the relative harmonic measures are the same from both sides. This means that for any measurable subset $E \subset e$,
$$
\lim_{z \to p, \, z \in \Omega_{p,X}} \frac{\omega_z(E)}{\omega_z(e)} = \lim_{z \to q, \, z \in \Omega_{q,Y}} \frac{\omega_z(E)}{\omega_z(e)}.
$$
\end{enumerate}

Arguing as in Section \ref{sec:convergence}, one can show that the true trees $\mathcal T_n$ converge to an infinite tree union the Julia set of $z^2+1/4$.

\subsection{Limit of Shabat polynomials}

We write $X$ for one of the symbols $L, R$.
We may further decompose each tile $\Omega_{p,X} \subset \Omega$ into countably many triangles $\triangle(e,p, X)$ by connecting the vertices in $\partial \Omega_{p,X}$ to the cusp $p \in \partial \Omega_{p,X}$ by hyperbolic geodesics in $\Omega_{p,X}$. 
We colour the triangles $\triangle(e, p, X) \subset \Omega$ black and white, so that 
$$
\triangle \, = \,
\triangle \bigl (\overline{v_{\Root}f(v_{\Root})}, 1/2, R \bigr )
\, \subset \, \Omega_{1/2,R} \, = \, \Omega_{1/2} \cap \mathbb{H}
$$ 
is white and adjacent triangles have different colours. Reflecting $\triangle$ in the real line, we get a triangle $\overline{\triangle} \subset
\Omega_{1/2,L} = \Omega_{1/2} \cap \mathbb{L}$. The union  $\triangle  \cup \overline{v_{\Root}f(v_{\Root})} \cup \overline{\triangle}$ constitutes a fundamental domain for the action of $f$ on $\Omega$.

\paragraph*{Mapping properties of the cosine.}
To describe the mapping properties of $\kappa(z) = \cos(\pi z)$, we draw the lines $\{y=0\}$ and $\{x = n: n \in \mathbb{Z} \}$ in the complex plane. These lines divide $\mathbb{C}$ into vertical half-strips $\{ \mathcal S_{n,\pm} \}$ of width 1. These may be coloured black and white so that adjacent half-strips have
opposite colours, with
$$
\mathcal S_{0,+} = \{ z \in \mathbb{C} \, : \, 0 <  \re z < 1, \, 0 < \im z < \infty \}
$$
being white. The map $\kappa$ takes each black half-strip conformally onto the upper half-plane and each white half-strip conformally onto the lower half-plane. The horizontal side of each $\mathcal S_{n, \pm}$ is mapped to the interval $[-1,1]$, while the vertical sides are mapped to the intervals $(-\infty, -1]$ and $[1,\infty)$.

\paragraph*{Mapping properties of the Fatou coordinate.}
The Fatou coordinate at the parabolic fixed point $1/2 \in \mathcal J(z^2+1/4)$ provides a conformal bijection between the quotient cylinder $\Omega/(z \sim f(z))$ and $\mathbb{C}/\mathbb{Z}$, which is  uniquely determined up to adding a constant in $\mathbb{C}/\mathbb{Z}$. Recall from the statement of Theorem \ref{main-thm2} that we use the normalization $\psi(v_{\Root}) = 0$.

\begin{lemma}
\label{psi} The Fatou coordinate $\psi$ is a holomorphic function on $\Omega$ which maps triangles $\triangle(e,p, X)$ conformally onto half-strips of the same colour.
\end{lemma}

\begin{proof}
Define $\psi_1: \triangle  \cup \overline{v_{\Root}f(v_{\Root})} \cup \overline{\triangle} \to \mathbb{C}$ to be conformal map which takes $\triangle$ to $\mathcal S_{0,+}$ and $\overline{\triangle}$ to $\mathcal S_{0,-}$
with
$$
v_{\Root} \to 0, \quad f(v_{\Root}) \to 1, \quad 1/2 \to \infty.
$$
The map $\psi_1$ extends to $\Omega$ using the functional equation $\psi_1(f(z)) = \psi_1(z) + 1$. Since $\psi_1$ possesses the properties that uniquely determine $\psi$, the two functions must be equal.
\end{proof}

Composing the above mappings, we get:

\begin{corollary}
\label{pi-psi}
 The map $z \to \cos(\pi \psi(z))$ takes each triangle $\triangle(e, p, X) \subset \Omega$ conformally onto the upper half-plane or the lower half-plane, with black triangles mapping onto the upper half-plane $\mathbb{H}$ and white triangles mapping onto the lower half-plane $\mathbb{L}$. Furthermore, $\cos(\pi \psi(z))$ takes edges to $[-1,1]$, cusps to infinity and $v_{\Root}$ to 1.
\end{corollary}

Considerations similar to the ones in Section \ref{sec:convergence-shabat} show that the the limit $h(z)$ of the Shabat polynomials $p_n(z)$ has the same description as the function $\cos(\pi \psi(z))$ described in Corollary \ref{pi-psi}. This completes the proof of Theorem \ref{main-thm2}.

\bibliographystyle{amsplain}

\newpage

 \begin{figure}[h!]
\centering
\includegraphics[scale=0.16]{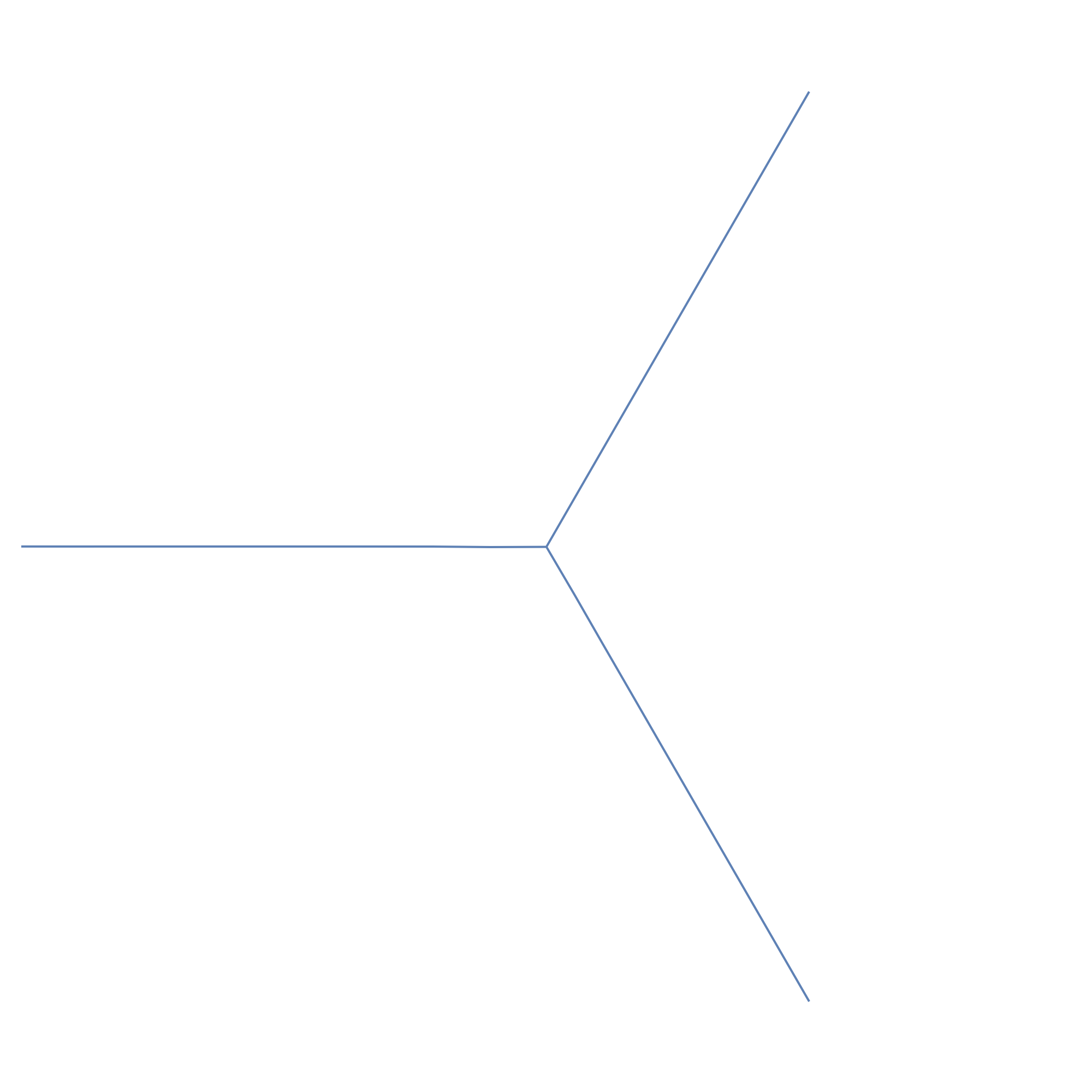}
\includegraphics[scale=0.16]{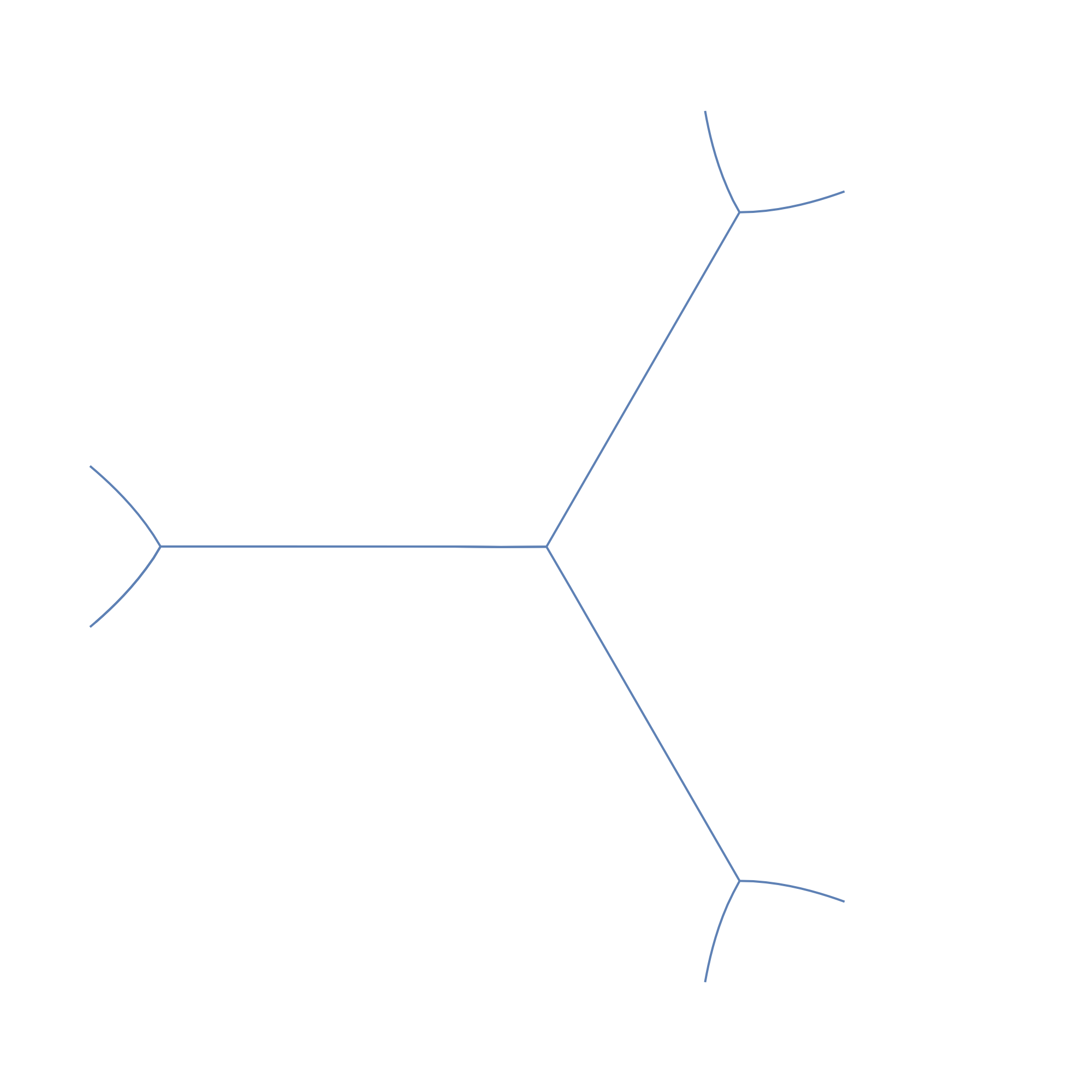}

\includegraphics[scale=0.16]{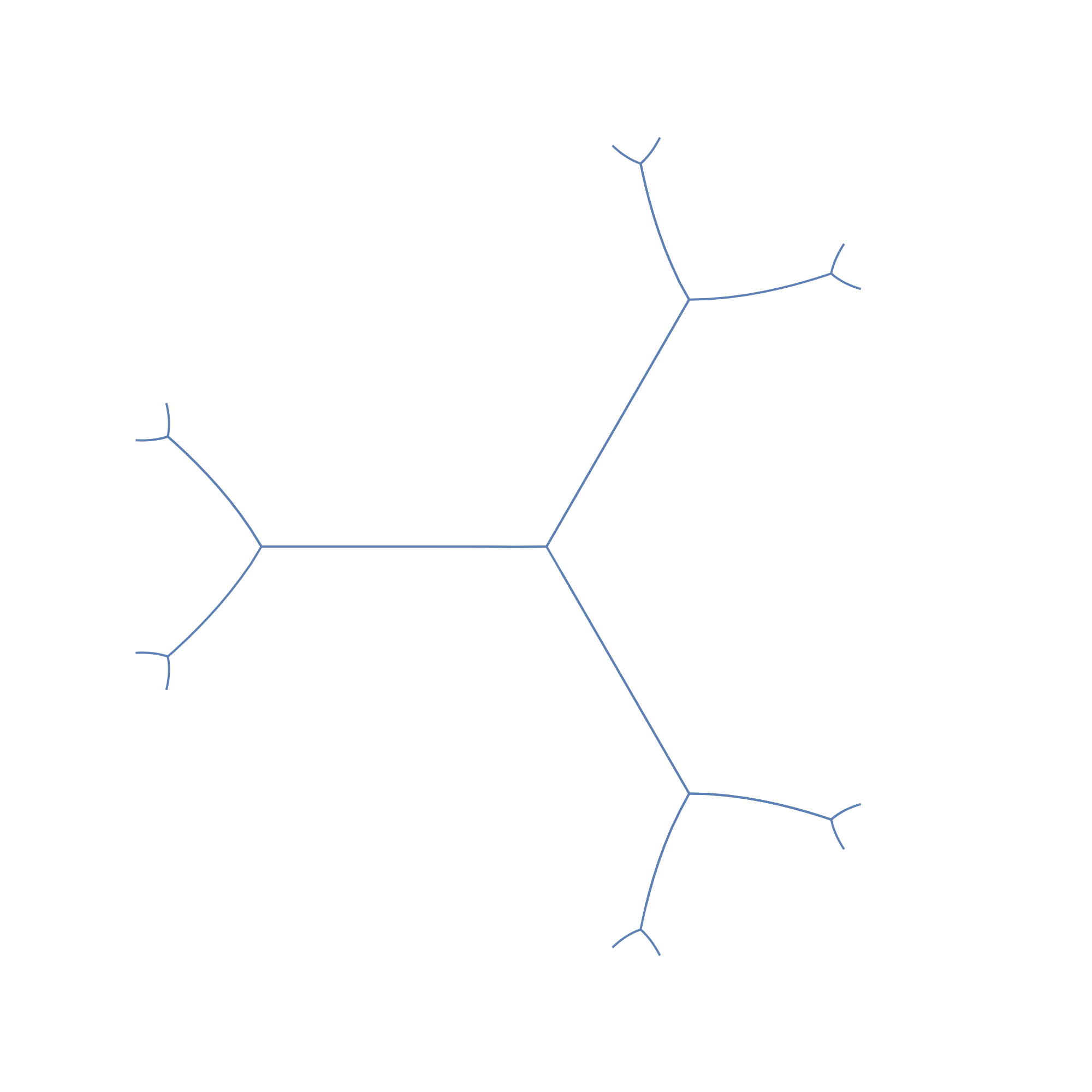}
\includegraphics[scale=0.16]{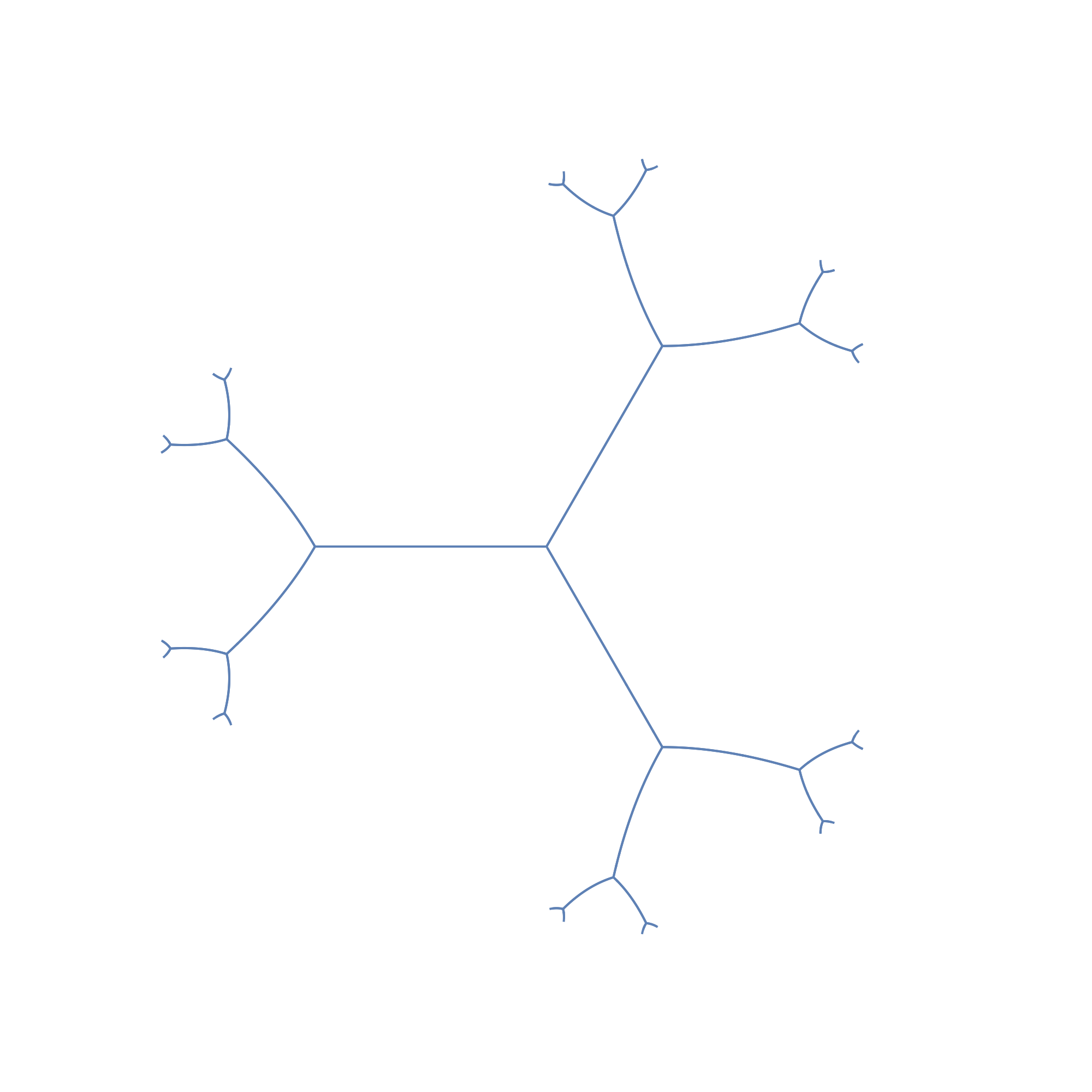}

\includegraphics[scale=0.16]{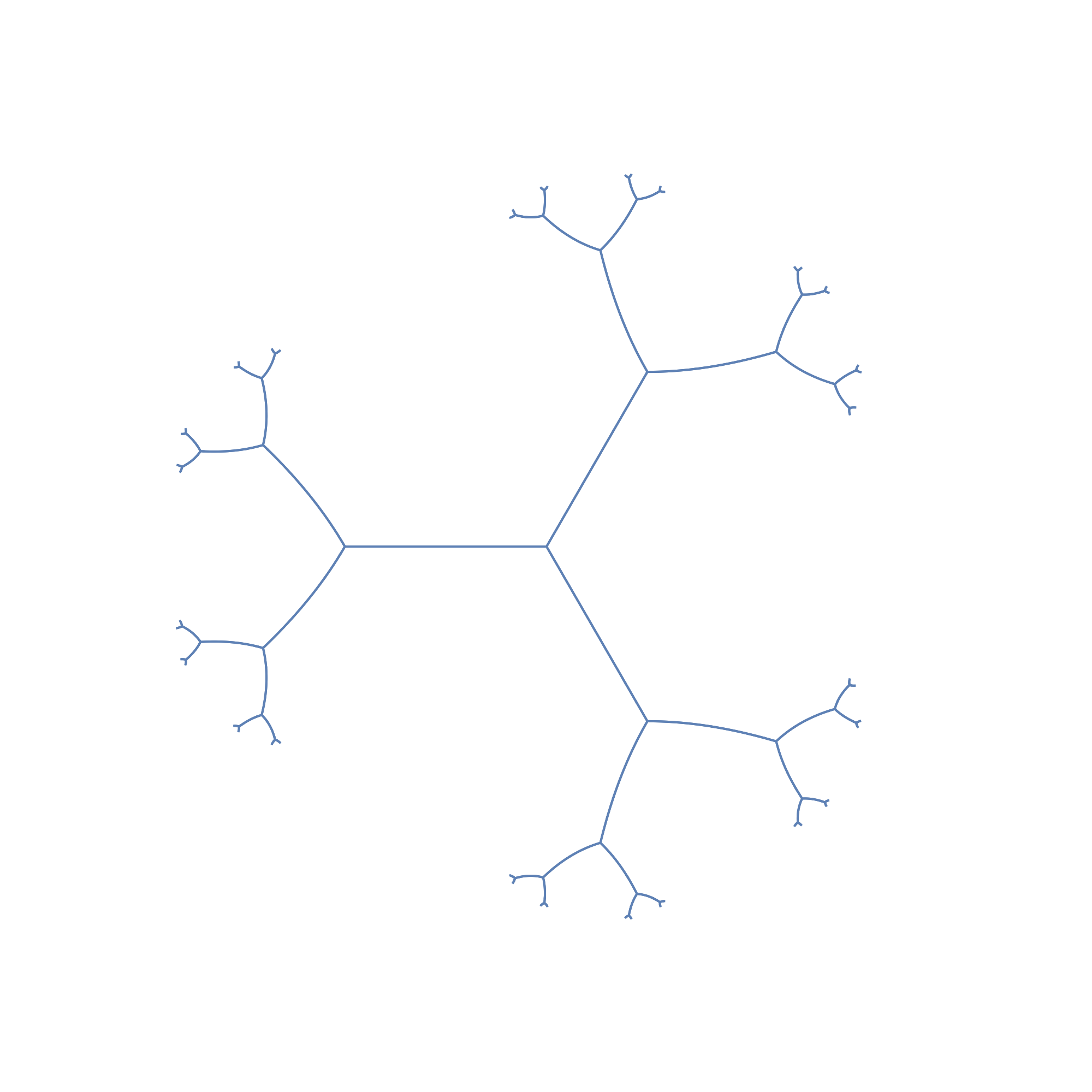}
\includegraphics[scale=0.16]{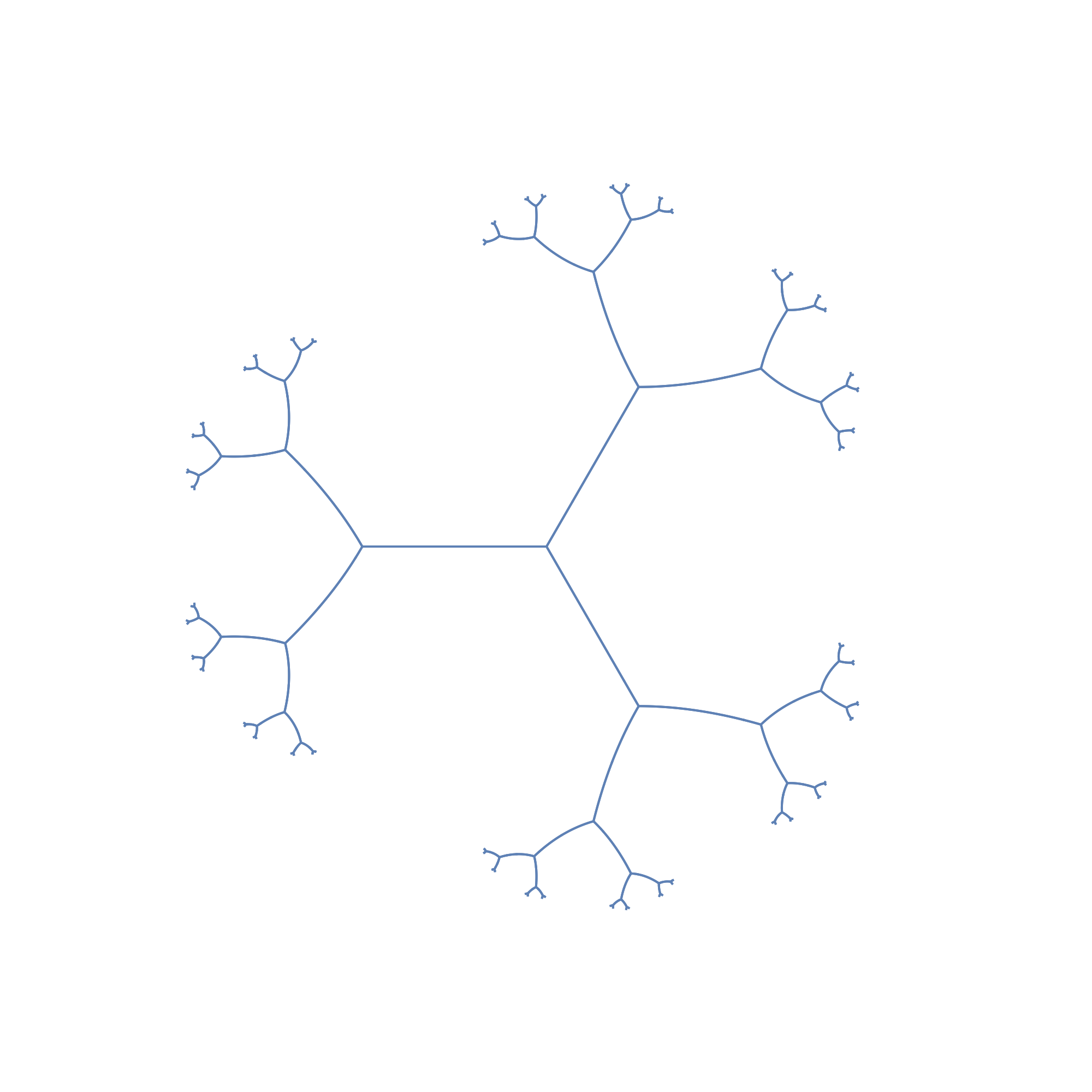}
\caption{A sequence of trees converging to the union of an infinite trivalent tree and the boundary of the developed deltoid.}
\label{fig:deltoid-trees}
 \end{figure}
 
 \newpage
 
  \begin{figure}[h!]
\centering
\includegraphics[scale=0.16]{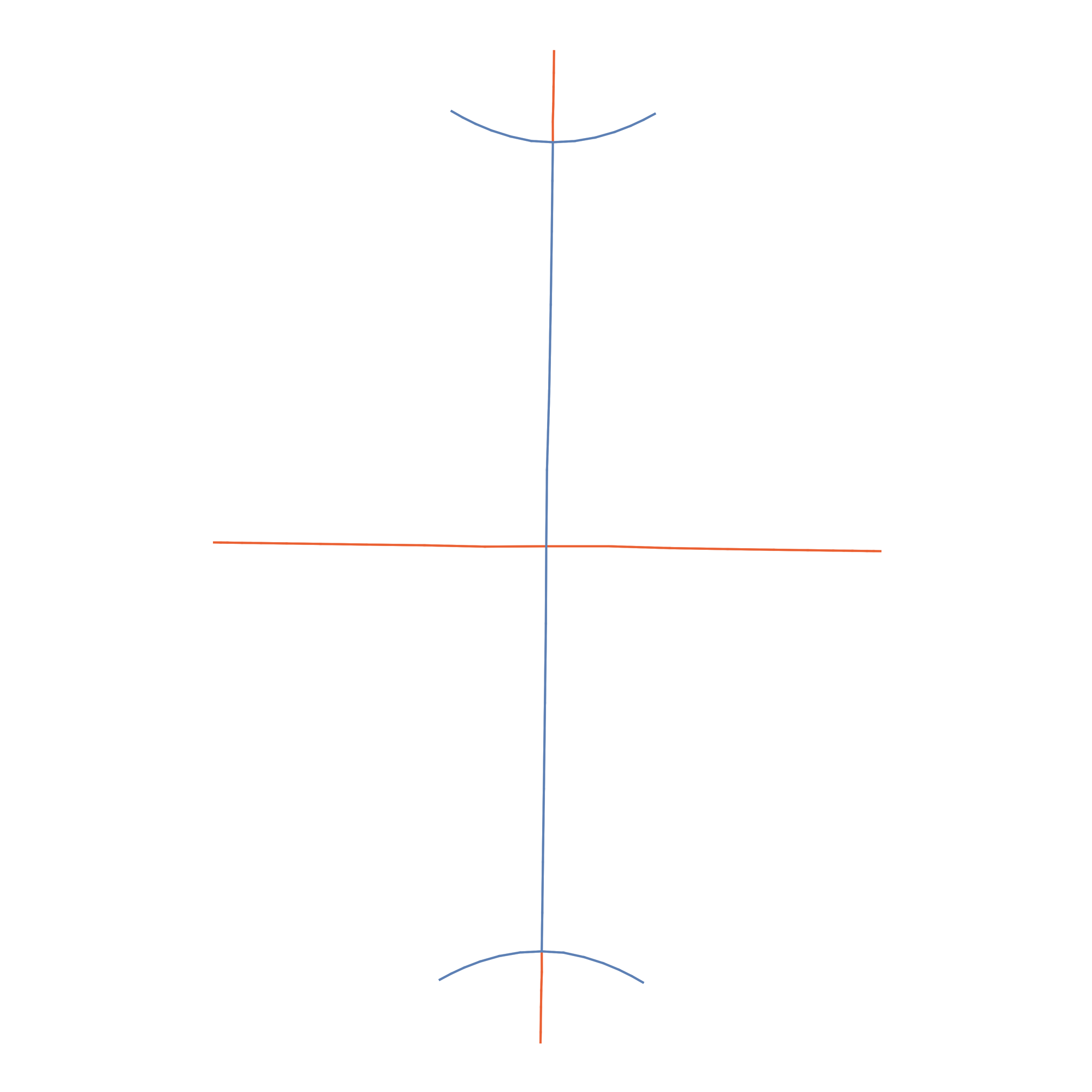}
\includegraphics[scale=0.16]{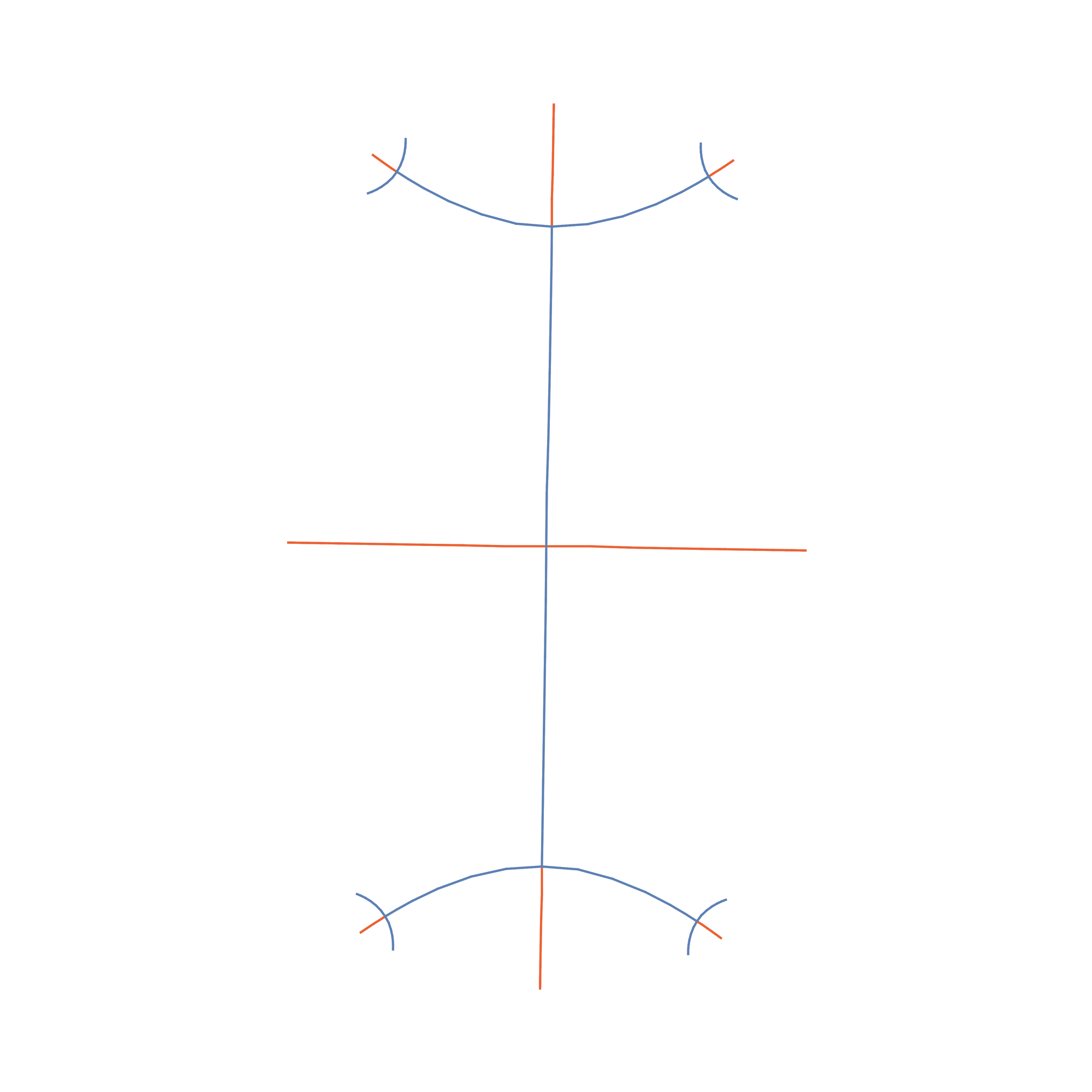}

\includegraphics[scale=0.16]{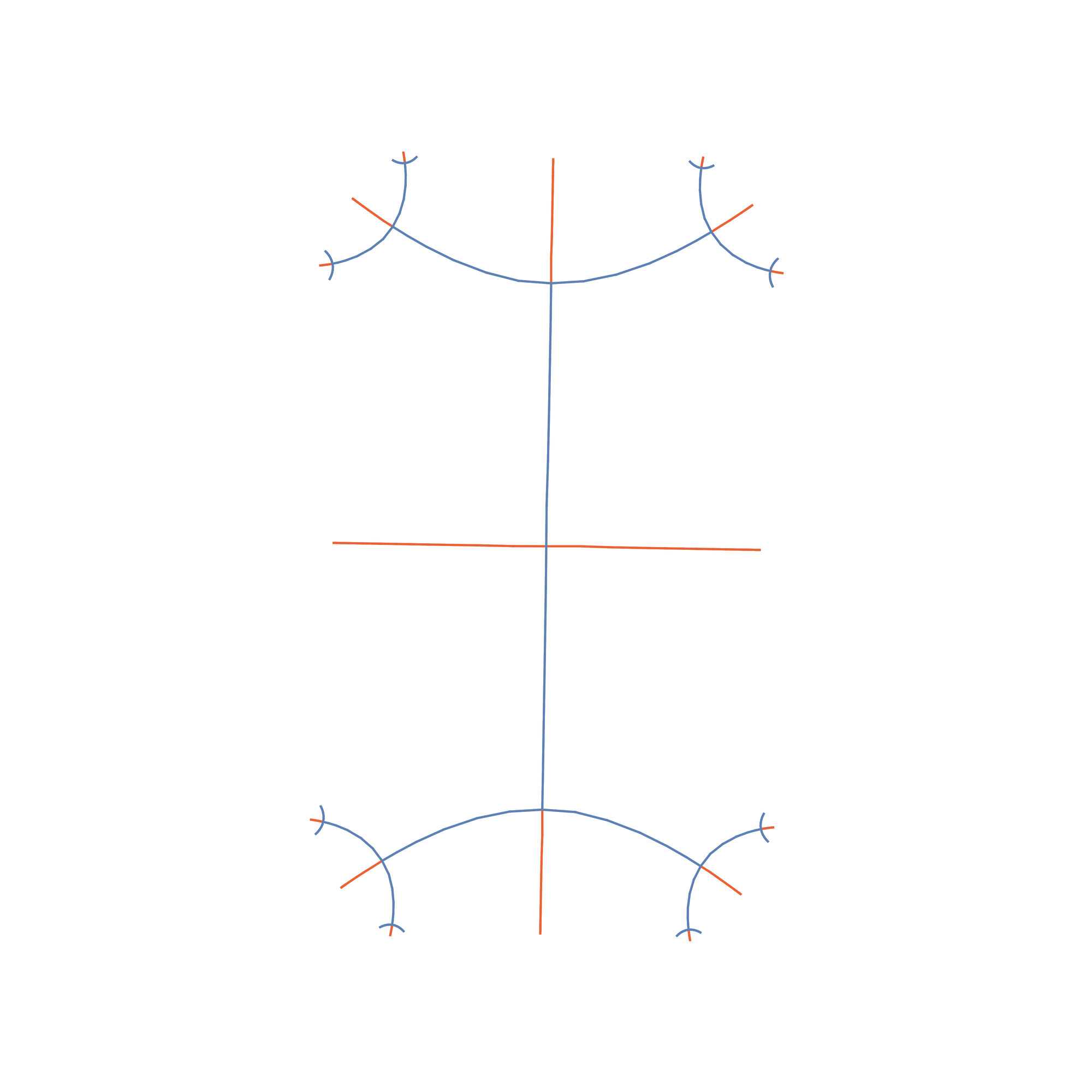}
\includegraphics[scale=0.16]{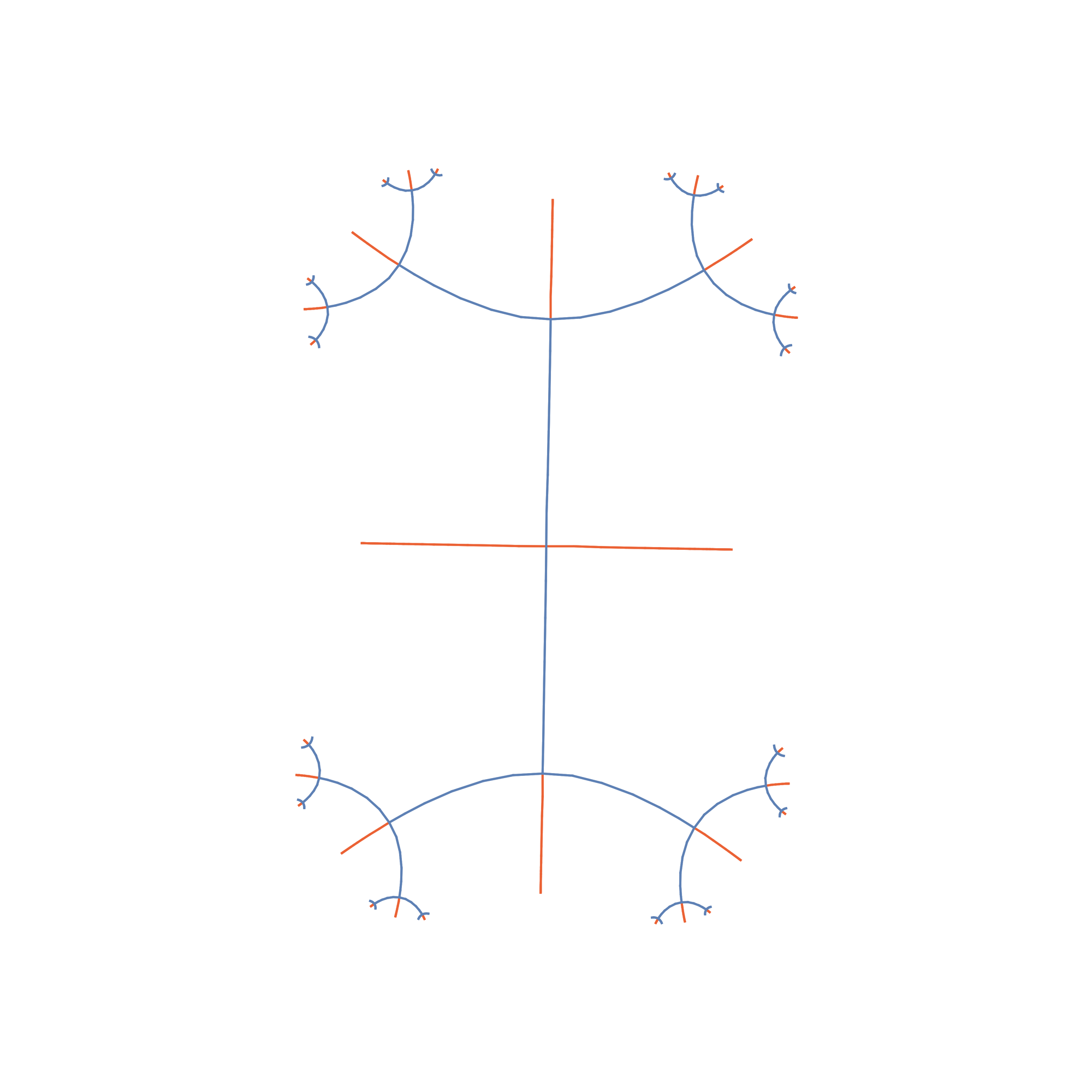}

\includegraphics[scale=0.16]{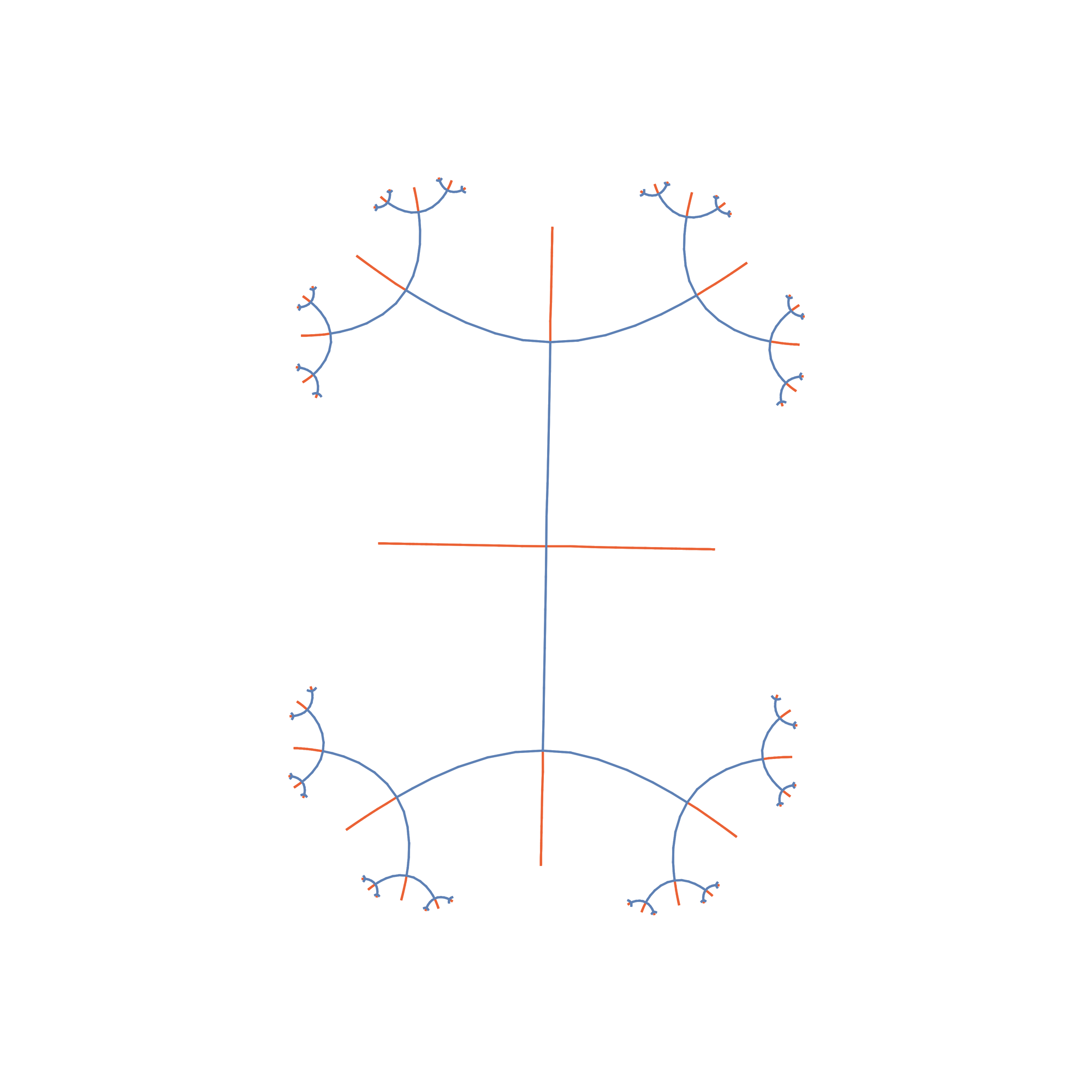}
\includegraphics[scale=0.16]{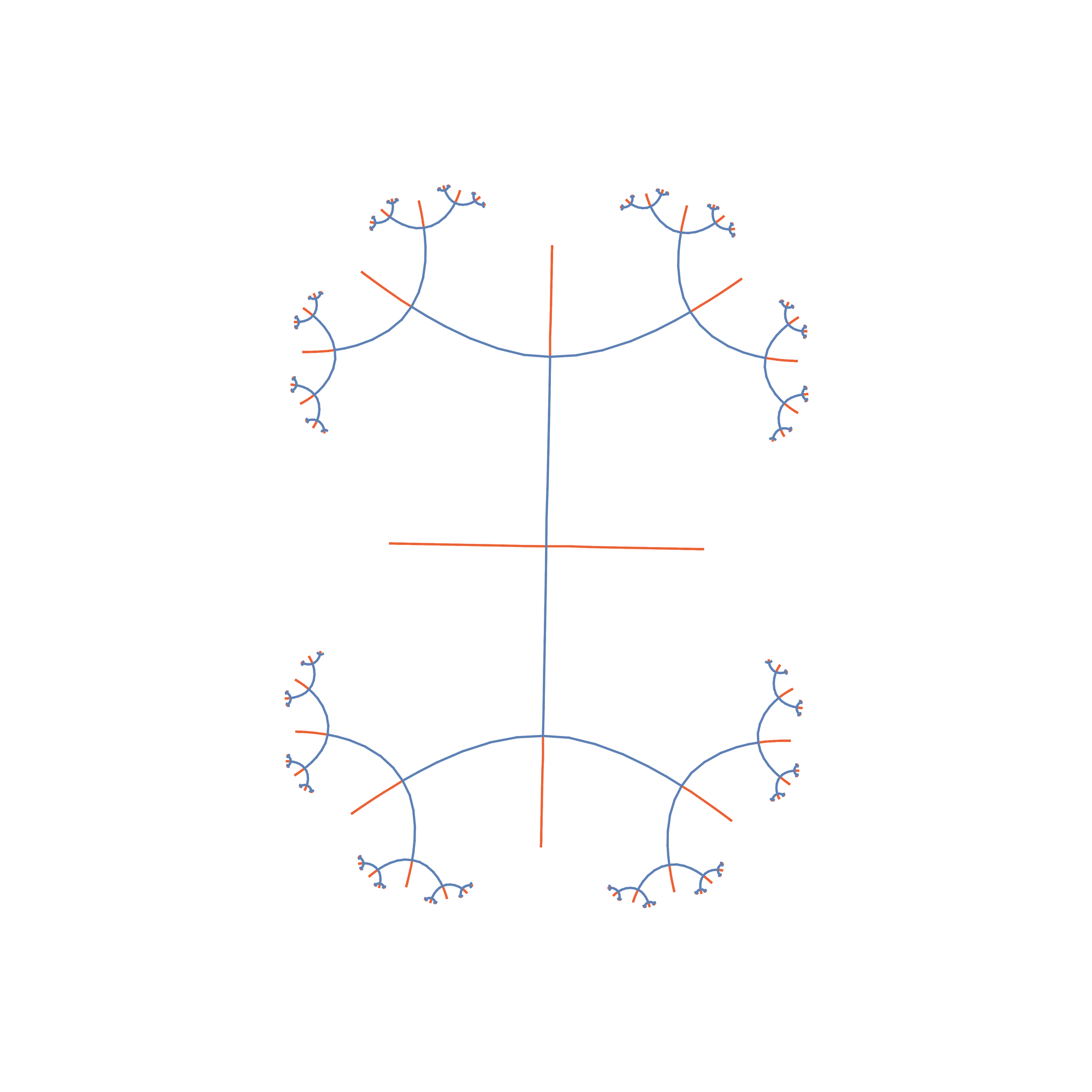}
\caption{A sequence of true trees converging to the union of the Julia set $\mathcal J(z^2+1/4)$ and an infinite tree whose vertex set is the grand orbit of the critical point $0$.}
\label{fig:cauliflower-trees}
 \end{figure}

\end{document}